\documentclass[11pt,a4paper]{article}
\title{\bf Generalized stochastic flow associated to the It\^o SDE with partially
Sobolev coefficients and applications}

\author{\sc{By Dejun Luo}
\vspace{3mm}\\
{\footnotesize Institute of Applied Mathematics, Chinese Academy of Sciences}}
\date{}
\usepackage{amssymb,amsmath,amsfonts,amsthm,array,color}

\def\R{\mathbb{R}}
\def\E{\mathbb{E}}
\def\P{\mathbb{P}}

\def\F{\mathcal{F}}
\def\L{\mathcal{L}}
\def\ul{\underline}
\def\ch{{\bf 1}}
\newcommand{\ra}{\rightarrow}

\newcommand{\da}{\downarrow}
\newcommand{\ua}{\uparrow}
\newcommand{\ee}{\varepsilon}

\def\div{\textup{div}}
\def\d{\textup{d}}
\def\supp{\textup{supp}}

\def\fin{\hfill$\square$}
\def\<{\langle}
\def\>{\rangle}
\def\ol{\overline}
\def\ul{\underline}

\def\Xint#1{\mathchoice
{\XXint\displaystyle\textstyle{#1}}%
{\XXint\textstyle\scriptstyle{#1}}%
{\XXint\scriptstyle\scriptscriptstyle{#1}}%
{\XXint\scriptscriptstyle\scriptscriptstyle{#1}}%
 \!\int}
\def\XXint#1#2#3{{\setbox0=\hbox{$#1{#2#3}{\int}$}
\vcenter{\hbox{$#2#3$}}\kern-.5\wd0}}

\def\dashint{\Xint-}
\def\bint{\dashint}

\begin{document}

\maketitle
\makeatletter 
\renewcommand\theequation{\thesection.\arabic{equation}}
\@addtoreset{equation}{section}
\makeatother 

\newtheoremstyle{newthm}
 {3pt}
 {3pt}
 {\itshape}
 {}
 {}
 {\textbf{.}}
 {.5em}
 {\thmname{\textbf{#1}}\thmnumber{ \textbf{#2}}\thmnote{ {\textbf{(#3)}}}}

\theoremstyle{newthm}

\newtheorem{theorem}{Theorem}[section]
\newtheorem{lemma}[theorem]{Lemma}       
\newtheorem{corollary}[theorem]{Corollary}
\newtheorem{proposition}[theorem]{Proposition}
\newtheorem{remark}[theorem]{Remark}
\newtheorem{example}[theorem]{Example}
\newtheorem{definition}[theorem]{Definition}

\begin{abstract}
We consider the It\^o SDEs with partially Sobolev coefficients.
Under some suitable conditions, we show the existence, uniqueness and stability
of generalized stochastic flows associated to such equations. As an application, we prove the weak
differentiability in the sense of measure of the stochastic flow generated by the It\^o SDE
with Sobolev coefficients.
\end{abstract}

\section{Introduction}

We consider the following stochastic differential equation
  \begin{equation}\label{Ito-SDE}
  \d X_t=\sigma(X_t)\,\d B_t+b(X_t)\,\d t,\quad X_0=x\in\R^n,
  \end{equation}
in which $\sigma=(\sigma^{ik})_{1\leq i\leq n,1\leq k\leq m}$ is a
matrix-valued function, $b=(b^1,\ldots,b^n)$ is a vector field,
and $B_t$ is an $m$-dimensional standard Brownian motion. It is well
known that if $\sigma$ and $b$ are globally Lipschitz continuous,
then equation \eqref{Ito-SDE} generates a unique stochastic flow of
homeomorphisms on $\R^n$. When the coefficients are less regular, for
instance, they only have log-Lipschitz continuity, it is still possible to prove
the existence of a homeomorphic flow, see \cite{Zhang05, FangLuo07}.

On the other hand, recently there are intensive studies on ODEs
  \begin{equation}\label{ODE}
  \frac{\d X_t}{\d t}=b(X_t),\quad X_0=x\in\R^n,
  \end{equation}
with weakly differentiable coefficients, see for instance
\cite{DiPernaLions89, Ambrosio04, CrippadeLellis}. Here by weakly differentiable
coefficients, we mean that they have Sobolev or even BV regularity.
The methods adopted in \cite{DiPernaLions89, Ambrosio04} are quite indirect,
in the sense that the authors first established the well-posedness
of the corresponding first order PDEs (transport equation or continuity
equation), from which they deduced the existence and uniqueness of
generalized flow of measurable maps associated to \eqref{ODE} (see also
\cite{CiprianoCruzeiro05} where the standard Gaussian measure $\gamma_n$
is taken as the reference measure). This strategy can
be seen as an extension of the classical characteristics method,
and is now widely called the DiPerna--Lions theory. In \cite{LeBrisLions04,
LeBrisLions08}, Le Bris and Lions made use of these ideas to study the
Fokker--Planck type equations with Sobolev coefficients; based on Ambrosio's
commutator estimate for BV vector fields, we slightly extend their results
to the case where the drift coefficient has only BV regularity, see \cite{Luo12}. The generalization
of this theory to the infinite dimensional Wiener space has been done in
\cite{AmbrosioFigalli09, FangLuo10}, see also \cite{Luo10} in which we
studied the Fokker--Planck type equations on the Wiener space. In \cite{Dumas},
the authors gave a rather sketchy argument of how to extend the DiPerna--Lions
theory to compact Riemannian manifolds; by proving a commutator estimate involving
the heat semi-group and Sobolev vector fields on manifolds, this theory was
recently generalized in \cite{FangLiLuo} to complete Riemannian manifolds under
suitable conditions on the lower bound of the Ricci curvature.  Using the
pointwise characterization of Sobolev functions, Crippa and de Lellis
gave in \cite{CrippadeLellis} direct proofs to many of the results in
the DiPerna--Lions theory.

It seems that DiPerna and Lions's original method does not work
for studying SDE \eqref{Ito-SDE}, as pointed out in the introduction of
\cite{Zhang11}. X. Zhang successfully implemented in \cite{Zhang10} the direct
method of Crippa and de Lellis to the It\^o SDE and proved the existence and
uniqueness of stochastic flow of maps generated by \eqref{Ito-SDE}.
A drawback of the main result in \cite[Theorem 2.6]{Zhang10} is the requirement that
$|\nabla \sigma|$ is bounded, a condition which is weakened
in \cite{Zhang12}. In \cite{FangLuoThalmaier} the authors took the standard
Gaussian measure $\gamma_n$ as the reference measure, and obtained similar
results under the exponential integrability of $|\nabla \sigma|^2,\,
|\div_{\gamma_n}(\sigma)|^2$ and $|\div_{\gamma_n}(b)|$. Here $\div_{\gamma_n}$
denotes the divergence with respect to the Gaussian measure $\gamma_n$.
Note that the exponential integrability of $|\nabla \sigma|^2$
is quite weak, but that of $|\div_{\gamma_n}(\sigma)|^2$ prevents us from
covering the classical case of globally Lipschitz coefficients, see
\cite[Theorem 1.2]{FangLuoThalmaier}. This is one of the reasons that we
do not take $\gamma_n$ as the reference measure in this paper. Another reason is that the
results in Lemma \ref{appendix-lem-2} do not hold for the Gaussian
measure $\gamma_n$. Here we also mention that we choose a finite measure
on $\R^n$ as the reference measure and assume the divergences of the coefficients
$\sigma$ and $b$ are exponentially integrable, hence they can be unbounded
(both locally and globally, see Theorem \ref{Extension} and \cite{FangLuoThalmaier, Zhang12}),
while the papers \cite{DiPernaLions89, Ambrosio04, CrippadeLellis} are set
in the framework of the Lebesgue measure, hence the authors naturally assume
that the divergence $\div(b)$ (or its negative part $[\div(b)]^-$) is
bounded.

The present work is motivated by \cite{LeBrisLions04, ALM, CrippadeLellis},
in which the authors studied the weak differentiability of the generalized
flow associated to the ODE \eqref{ODE} with Sobolev vector field $b$. Again
the results in \cite{LeBrisLions04} are derived from the related transport
equation, while the ones in \cite{ALM, CrippadeLellis} follow from the
pointwise inequality of Sobolev functions. Since the generalized stochastic
flow of measurable maps has already been established in \cite{Zhang10,
FangLuoThalmaier, Zhang12}, we intend to study in this work the differentiability of
the stochastic flow. However, we are unable to transfer the methods in
\cite{ALM, CrippadeLellis} to the case of SDE for proving the approximate differentiability
of the stochastic flow. The main problem is that the level set $G_R$
(see Lemma \ref{sect-2-lem-1}) of the stochastic flow depends on the random
element $\omega$, hence one has to take expectation twice in order to estimate
an quantity of the form (2.5) in \cite{CrippadeLellis}. We do not know
how to handle this problem.

Therefore, we follow the idea of \cite{LeBrisLions04} to study the differentiability
in the sense of measure of the stochastic flow. To this end, we first consider a
special form of SDE \eqref{Ito-SDE} whose coefficients $\sigma$ and $b$ have the
structure below: there is $n_1\in\{1,\ldots,n-1\}$, such that
  $$\sigma_1:=(\sigma^{ij})_{1\leq i\leq n_1,1\leq j\leq m}
  \quad \mbox{and} \quad b_1:=(b^1,\ldots,b^{n_1})$$
only depend on the first $n_1$-variables $(x^1,\ldots,x^{n_1})$. In the following
we also denote by $\sigma_2$ (resp. $b_2$) the last $(n-n_1)$-rows (resp. components)
of the diffusion matrix $\sigma$ (resp. the drift $b$), and
$x_1=(x^1,\ldots,x^{n_1}),\,x_2=(x^{n_1+1},\ldots,x^n)$ (thus $x\in\R^n$ can be
written as $(x_1,x_2)$). Our basic assumptions,
among other conditions that will be specified later, are
  \begin{equation}\label{Regularity-1}
  \sigma_1\in W^{1,2q}_{x_1,loc},\quad b_1\in W^{1,q}_{x_1,loc};
  \end{equation}
and
  \begin{equation}\label{Regularity-2}
  \sigma_2\in L^{2q}_{x_1,loc}(W^{1,2q}_{x_2,loc}),\quad b_2\in L^q_{x_1,loc}(W^{1,q}_{x_2,loc}).
  \end{equation}
Here $q>1$ is a fixed number. Note that we don't require $\sigma_2$ and $b_2$ have Sobolev
regularity with respect to $x_1$.

The paper is organized as follows. In Section 2 we first recall the definition of
generalized stochastic flow associated to It\^o's SDE \eqref{Ito-SDE}.
After that, we extend the known results
on the existence and uniqueness of stochastic flows generated by It\^o's SDE
to allow the coefficients to be locally unbounded. Recall that the main results
in \cite{Zhang10, FangLuoThalmaier, Zhang12} require the coefficients $\sigma$
and $b$ have linear growth. This extension is necessary for proving the
differentiability of the stochastic flow, since the linear growth condition
for the second equation in \eqref{sect-5.2} will basically result in
the boundedness of the gradients of $\sigma$ and $b$, which is too restrictive.

Then we state and prove an intermediate result in Section 3, where the coefficients
$\sigma_2\in W^{1,2q}_{x_1,x_2,loc}$ and $b_2\in W^{1,q}_{x_1,x_2,loc}$. One reason for
establishing such a result is to avoid the regularization of the coefficients
$\sigma_1$ and $b_1$ in the proof of the existence of stochastic flows generated by
It\^o's SDE with partially Sobolev coefficients (see Theorem \ref{sect-4-existence});
otherwise, we cannot apply the a-priori estimate in Lemma \ref{a-priori-estimate},
since the coefficients $\sigma_2$ and $b_2$ have no Sobolev regularity on
the variable $x_1=(x^1,\ldots,x^{n_1})$. We also
find a uniform estimate of the Radon--Nikodym
density of the form Lemma \ref{sect-3-lem-3}, which does not involve the
exponential integrability of $|\nabla_{x_1}\sigma_2|^2$.

The main result of this paper is presented in Section 4, in which
the key step is to prove an a-priori estimate which follows the idea of Crippa and
de Lellis \cite[Theorem 3.8]{CrippadeLellis} and has appeared in \cite{Zhang10,
FangLuoThalmaier, Zhang12} in similar forms. The main difference between this
estimate and the previous ones is that we only assume partial Sobolev
regularity on the coefficients. As some of the arguments
in Sections 3 and 4 are analogous to those of Section 2, we only give
relatively detailed proofs in Section 2 and omit them in the subsequent
sections to save space.

In Section 5 we apply the results obtained in the previous section to show the
weak differentiability in the sense of measure of the generalized stochastic
flow of measurable maps, following the ideas in \cite[Section 4]{LeBrisLions04}. The main
part consists in checking that the systems of It\^o equations fulfil the assumptions in
Section 4.

Finally, we present in the appendix some preliminary
results that are frequently used in the paper. Especially, we
give a careful analysis of the expression of the Radon--Nikodym density which
makes it possible for us to study the SDE with the above-mentioned
special structure. We also prove an inequality for the integral of
local maximal functions on the whole $\R^n$ with respect to some general
finite measure which seems to have independent interest.

\section{The It\^o SDE with locally unbounded coefficients}

First of all we give the precise meaning of the generalized stochastic flow
(cf. \cite[Definition 5.1]{FangLuoThalmaier} and \cite[Definition 2.1]{Zhang12}).
This notion is related to some reference measure on $\R^n$. In this paper,
we mainly consider the following type of measures: for some $\alpha>n/2$, set
  \begin{equation}\label{reference-measure}
  \lambda(x)=-\alpha\log(1+|x|^2)\ (x\in\R^n)\quad
  \mbox{and}\quad  \d\mu=e^{\lambda(x)}\,\d x.
  \end{equation}
The exact value of $\alpha$ has no importance. It is clear that $\mu(\R^n)<+\infty$.
Denote by $\theta_sB$ the time-shift of the Brownian motion, that is,
$(\theta_sB)_t=B_{t+s}-B_s$ for all $t\geq0$. For a measurable map
$\varphi:\R^n\ra\R^n$, we write $\varphi_\#\mu=\mu\circ\varphi^{-1}$
for the push-forward of $\mu$ by $\varphi$ (also called the distribution
of $\varphi$ under $\mu$).

\begin{definition}\label{sect-2-def}
We say that a measurable map $X\colon\Omega\times\R^d\ra C([0,T],\R^n)$
is a generalized stochastic flow associated to the It\^{o} SDE \eqref{Ito-SDE} if
\begin{enumerate}
\item[\rm(i)] for each $t\in [0,T]$ and almost all $x\in\R^n$, $\omega\ra
X_t(\omega,x)$ is measurable with respect to $\F_t$, i.e., the natural
filtration generated by the Brownian motion $\{B_s\colon s\leq t\}$;

\item[\rm(ii)] for each $t\in [0,T]$, there exists $K_t\in
L^1(\P\times\mu)$ such that $(X_t(\omega,\cdot))_\#\mu$ admits
$K_t$ as the density with respect to $\mu$;

\item[\rm(iii)] for $(\P\times\mu)$-a.e. $(\omega,x)$,
  \begin{equation*}
  \int_0^T|\sigma(X_s(\omega,x))|^2\,\d s+\int_0^T|b(X_s(\omega,x))|\,\d
  s<+\infty;
  \end{equation*}
\item[\rm(iv)] for $\mu$-a.e. $x\in\R^n$, the integral equation below holds
almost surely:
  \begin{equation*}
  X_t(\omega,x)=x+\int_0^t \sigma(X_s(\omega,x))\,\d B_s+\int_0^t
  b(X_s(\omega,x))\,\d s,\quad \mbox{for all } t\in [0,T];
  \end{equation*}
\item[\rm(v)] the flow property holds
  \begin{equation*}
  X_{t+s}(\omega,x)=X_t(\theta_sB, X_s(\omega,x)).
  \end{equation*}
\end{enumerate}
\end{definition}

In this section we slightly extend the main results of \cite{Zhang10, FangLuoThalmaier,
Zhang12} to allow the coefficients $\sigma$ and $b$ to be locally unbounded, while
the aforementioned papers required that the coefficients have linear growth. To this end,
we introduce some notations. Fix some $q>1$ and take $\alpha>q+n/2$ in the definition
\eqref{reference-measure} of the reference measure. We also denote by $\bar\sigma=\frac{\sigma}{1+|x|}$
and $\bar b=\frac{b}{1+|x|}$ to simplify the notations. We assume the following conditions:
  \begin{itemize}
  \item[\rm(C1)] $\sigma\in W^{1,2q}_{loc},b\in W^{1,q}_{loc}$;
  \item[\rm(C2)] there is a $p_0>0$ such that
  $\int_{\R^{n}}\exp\big[p_0\big([\div(b)]^- +|\bar b|+|\bar\sigma|^2
  +|\nabla\sigma|^2\big)\big]\d\mu<+\infty$.
  \end{itemize}

\begin{remark}\label{sect-2-rem} We have the following observations.
\begin{itemize}
\item[\rm(i)] It is clear that when $\sigma$ and $b$ are globally Lipschitz
continuous, they satisfy the conditions (C1) and (C2).
\item[\rm(ii)] The condition (C2) implies $\bar\sigma,\bar b\in L^p(\mu)$
for any $p>1$. By the choice of $\alpha$, there is $p$ sufficiently big such that
$2\alpha-n>2qp/(p-1)$, hence $\int_{\R^n}(1+|x|)^{2qp/(p-1)}\d\mu<+\infty$. By
H\"older's inequality,
  $$\int_{\R^n}|\sigma|^{2q}\,\d\mu\leq \bigg[\int_{\R^n}|\bar\sigma|^{2qp}\,\d\mu\bigg]^{1/p}
  \bigg[\int_{\R^n}(1+|x|)^{2qp/(p-1)}\d\mu\bigg]^{(p-1)/p}<+\infty.$$
Thus $\sigma\in L^{2q}(\mu)$. In the same way we have $b\in L^{2q}(\mu)$.
\item[\rm(iii)] Suppose that $\supp(b)\subset B(1)$ and there is $\beta\in (0,n/p_0)$
such that $|\bar b(x)|\leq \log\frac1{|x|^\beta}$
for all $|x|\leq 1$, then $\int_{\R^n}e^{p_0|\bar b|}\,\d\mu<+\infty$.
Hence the coefficient $b$ (and also $\sigma$) of the It\^o SDE can be locally unbounded.
\end{itemize}
\end{remark}

We shall prove

\begin{theorem}\label{Extension} Under the conditions (C1) and (C2), there exists
a unique generalized stochastic flow associated to the It\^o SDE \eqref{Ito-SDE}.
Moreover, the Radon--Nikodym density $\rho_t$ of the flow with respect to the
reference measure $\mu$ satisfies $\rho_t\in L^1\log L^1$.
\end{theorem}

Here by $\rho_t\in L^1\log L^1$ we mean that
$\E\int_{\R^n}\rho_t|\log\rho_t|\,\d\mu<+\infty$.
We remark that when $t$ is small enough, the flow $X_t$ is integrable on $\R^n$ with
respect to $\mu$, which is an easy consequence of Lemma \ref{sect-2-lem-1}
and Proposition \ref{sect-2-prop-2}. The integrability of $X_t$ for general $t>0$ can be proved
if we strengthen the condition (C2) by requiring that it holds for any $p_0>0$;
however, this condition is too restrictive.

We shall divide the proof of this theorem into several steps, which are presented
in the following lemmas and propositions. First we
prove an a-priori estimate on the level set of the solution flow $X_t$.
We denote by $\|\cdot\|_{\infty,T}$ the supremum norm in $C([0,T],\R^n)$, the
space of continuous curves in $\R^n$. For $R>0$, define the level set
  $$G_R=\big\{(\omega,x)\in\Omega\times\R^n:\|X_\cdot(\omega,x)\|_{\infty,T}\leq R\big\}.$$

\begin{lemma}[Estimate of level sets]\label{sect-2-lem-1}
Let $X_t$ be a generalized stochastic flow associated to It\^o SDE \eqref{Ito-SDE},
and $\rho_t$ the Radon--Nikodym density with respect to $\mu$. Suppose that
  $$\Lambda_{p,T}:=\sup_{0\leq t\leq T}\|\rho_t\|_{L^p(\P\times\mu)}<+\infty,$$
where $p$ is the conjugate number of $q$. Then under the condition (C2), we have
  $$(\P\times\mu)(G_R^c)\leq \frac{C}{R},$$
where $C$ depends on $T,\Lambda_{p,T}$, $\|\sigma\|_{L^{2q}(\mu)}$
and $\|b\|_{L^{q}(\mu)}$.
\end{lemma}

\noindent{\bf Proof.} First we deduce from (C2) and Remark \ref{sect-2-rem}(ii)
that $\|\sigma\|_{L^{2q}(\mu)}$ and $\|b\|_{L^{q}(\mu)}$ are finite.
For a.e. $(\omega,x)\in \Omega\times\R^n$, we have
  $$X_t(x)=x+\int_0^t\sigma(X_s(x))\,\d B_s+\int_0^t b(X_s(x))\,\d s.$$
Therefore
  \begin{equation}\label{sect-2-lem-1-1}
  \|X_\cdot(x)\|_{\infty,T}\leq |x|+\sup_{0\leq t\leq T}\bigg|\int_0^t\sigma(X_s(x))\,\d B_s\bigg|
  +\sup_{0\leq t\leq T}\bigg|\int_0^t b(X_s(x))\,\d s\bigg|.
  \end{equation}
By Burkholder's inequality,
  \begin{align*}
  \E\sup_{0\leq t\leq T}\bigg|\int_0^t\sigma(X_s(x))\,\d B_s\bigg|
  \leq 2\bigg[\E\int_0^T|\sigma(X_s(x))|^2\,\d s\bigg]^{\frac12}.
  \end{align*}
Now Cauchy's inequality leads to
  \begin{align*}
  \int_{\R^n}\E\sup_{0\leq t\leq T}\bigg|\int_0^t\sigma(X_s(x))\,\d B_s\bigg|\d\mu
  &\leq 2\mu(\R^n)^{\frac12}\bigg[\int_0^T\E\int_{\R^n}|\sigma(X_s(x))|^2\,\d\mu(x)\d s\bigg]^{\frac12}\cr
  &= 2\mu(\R^n)^{\frac12}\bigg[\int_0^T\E\int_{\R^n}|\sigma(y)|^2\rho_s(y)\,\d\mu(y)\d s\bigg]^{\frac12}.
  \end{align*}
We have by H\"older's inequality that
  $$\E\int_{\R^n}|\sigma(y)|^2\rho_s(y)\,\d\mu(y)\leq
  \|\sigma\|_{L^{2q}(\mu)}^2\|\rho_s\|_{L^p(\P\times\mu)}
  \leq \Lambda_{p,T}\|\sigma\|_{L^{2q}(\mu)}^2.$$
Therefore
  \begin{equation}\label{sect-2-lem-1-2}
  \int_{\R^n}\E\sup_{0\leq t\leq T}\bigg|\int_0^t\sigma(X_s(x))\,\d B_s\bigg|\d\mu
  \leq 2(\mu(\R^n)\, T\Lambda_{p,T})^{\frac12}\|\sigma\|_{L^{2q}(\mu)}.
  \end{equation}

Next
  \begin{align*}
  \E\int_{\R^n}\sup_{0\leq t\leq T}\bigg|\int_0^t b(X_s(x))\,\d s\bigg|\d\mu
  &\leq \int_0^T\E\int_{\R^n}|b(X_s(x))|\,\d\mu(x)\d s\cr
  &= \int_0^T\E\int_{\R^n}|b(y)|\rho_s(y)\,\d\mu(y)\d s.
  \end{align*}
Again by H\"older's inequality,
  \begin{equation}\label{sect-2-lem-1-3}
  \E\int_{\R^n}\sup_{0\leq t\leq T}\bigg|\int_0^t b(X_s(x))\,\d s\bigg|\d\mu
  \leq \int_0^T\|b\|_{L^q(\mu)}\|\rho_s\|_{L^p(\P\times\mu)}\,\d s
  \leq T\Lambda_{p,T}\|b\|_{L^q(\mu)}.
  \end{equation}
Now integrating both sides of \eqref{sect-2-lem-1-1} on $\Omega\times\R^n$ and by
\eqref{sect-2-lem-1-2}, \eqref{sect-2-lem-1-3}, we get
  \begin{equation}\label{sect-2-lem-1-4}
  \E\int_{\R^n}\|X_\cdot(x)\|_{\infty,T}\,\d\mu\leq C_1
  +2(\mu(\R^n)\, T\Lambda_{p,T})^{\frac12}\|\sigma\|_{L^{2q}(\mu)}+T\Lambda_{p,T}\|b\|_{L^q(\mu)},
  \end{equation}
where $C_1:=\int_{\R^n}|x|\,\d\mu(x)<+\infty$. Finally by Chebyshev's inequality,
  $$(\P\times\mu)(G_R^c)\leq \frac1R \int_{\Omega\times\R^n}
  \|X_\cdot(x)\|_{\infty,T}\,\d(\P\times\mu)\leq \frac CR,$$
where $C$ is given by the right hand side of \eqref{sect-2-lem-1-4}. \fin

\medskip

Similar to \cite[Lemma 6.1]{Zhang10}, \cite[Theorem 5.2]{FangLuoThalmaier}
and \cite[Lemma 4.1]{Zhang12}, we have the following

\begin{lemma}[Stability estimate]\label{sect-2-lem-2}
Suppose that $\sigma,\tilde\sigma\in W^{1,2q}_{loc}$ and $b,\tilde b\in W^{1,q}_{loc}$. Let $X_t$
(resp. $\tilde X_t$) be the stochastic flow associated to the It\^o SDE \eqref{Ito-SDE}
with coefficients $\sigma$ and $b$ (resp. $\tilde\sigma$ and $\tilde b$). Denote by
$\rho_t$ (resp. $\tilde \rho_t$) the Radon--Nikodym density of $X_t$ (resp. $\tilde X_t$)
with respect to $\mu$. Assume that
  $$\Lambda_{p,T}:=\sup_{0\leq t\leq T}\big(\|\rho_t\|_{L^p(\P\times\mu)}\vee
  \|\tilde \rho_t\|_{L^p(\P\times\mu)}\big)<+\infty.$$
where $p$ is the conjugate number of $q$. Then for any $\delta>0$,
  \begin{align*}
  &\E\int_{G_R\cap \tilde G_R}\log\bigg(\frac{\|X-\tilde X\|^2_{\infty,T}}{\delta^2}+1\bigg)\d\mu\cr
  &\hskip6pt \leq C_T\Lambda_{p,T}\bigg\{C_{n,q}\Big[\|\nabla b\|_{L^q(B(3R))}
  +\|\nabla \sigma\|_{L^{2q}(B(3R))}+\|\nabla\sigma\|^2_{L^{2q}(B(3R))}\Big]\cr
  &\hskip60pt +\frac1{\delta^2}\|\sigma-\tilde\sigma\|^2_{L^{2q}(B(R))}
  +\frac1\delta \Big[\|\sigma-\tilde\sigma\|_{L^{2q}(B(R))}+
  \|b-\tilde b\|_{L^{q}(B(R))}\Big]\bigg\},
  \end{align*}
where $\tilde G_R:=\big\{(\omega,x)\in \Omega\times\R^n:\|\tilde X_\cdot(\omega,x)\|_{\infty,T}\leq R\big\}$
is the level set of the flow $\tilde X_t$.
\end{lemma}

Here the space $L^q(B(R))$ is defined with respect to the Lebesgue measure.
The proof of Lemma \ref{sect-2-lem-2} is similar to the above cited references, hence we omit it.

Now we start to prove the existence part of Theorem \ref{Extension}.
We have to regularize the coefficients $\sigma$ and $b$. Let
$\chi\in C_c^\infty(\R^n,\R_+)$ be such that $\int_{\R^n}\chi\,\d x=1$
and its support $\supp(\chi)\subset B(1)$. For $k\geq 1$, define $\chi_k(x)=k^n\chi(kx)$ for all
$x\in\R^n$. Next choose $\psi\in C_c^\infty(\R^n,[0,1])$ which satisfies
$\psi|_{B(1)}\equiv 1$ and $\supp(\psi)\subset B(2)$. Set $\psi_k(x)=\psi(x/k)$ for
all $x\in\R^n$ and $k\geq1$. Now we define
  $$\sigma_k=(\sigma\ast\chi_k)\,\psi_k
  \quad \mbox{and} \quad b_k=(b\ast\chi_k)\,\psi_k.$$
Then for every $k\geq1$, the functions $\sigma_k$ and $b_k$ are smooth with
compact supports. Consider the following It\^o's SDE:
  \begin{equation}\label{smooth-SDE}
  \d X^k_t=\sigma_k(X^k_t)\,\d B_t+b_k(X^k_t)\,\d t,\quad X^k_0=x.
  \end{equation}
This equation has a unique strong solution which gives rise to a
stochastic flow of diffeomorphisms on $\R^n$. Denote by $\rho^k_t$ the
Radon--Nikodym density of $(X^k_t)_\#\mu$ with respect to $\mu$. Applying
Lemma \ref{density-estimate} for $p>1$, we have
  \begin{equation}\label{sect-2.2}
  \|\rho^k_t\|_{L^p(\P\times\mu)}\leq \mu(\R^n)^{\frac1{p+1}}
  \bigg(\sup_{t\in[0,T]}\int_{\R^n}\exp\big(p^3t|\Lambda_1^{\sigma_k}|^2
  -p^2t\Lambda_2^{\sigma_k,b_k}\big)\d\mu\bigg)^{\frac1{p(p+1)}}.
  \end{equation}
We shall give a uniform estimate to the density functions. For this purpose we need

\begin{lemma}\label{sect-2-lem-3} There is a constant $C_0>0$, independent of
$k\geq 1$, such that
  \begin{itemize}
  \item[\rm(1)] $|\Lambda_1^{\sigma_k}|^2\leq C_0\big(|\div(\sigma)|^2
  +|\bar\sigma|^2\big)\ast\chi_k$;
  \item[\rm(2)] $-\Lambda_2^{\sigma_k,b_k}\leq C_0\big([\div(b)]^- +|\bar b|
  +|\nabla\sigma|^2+|\bar\sigma|^2\big)\ast\chi_k$.
  \end{itemize}
\end{lemma}

\noindent{\bf Proof.} (1) By the definition of $\Lambda_1^{\sigma_k}$, we have
  \begin{equation}\label{sect-2-lem-3.1}
  \Lambda_1^{\sigma_k}=\div(\sigma_k)+\sigma_k^\ast\nabla\lambda,
  \end{equation}
where $\sigma_k^\ast$ is the transpose of $\sigma_k$. For every $l\in\{1,\ldots,m\}$, we have
  $$\div(\sigma_k^{\cdot,l})=[\div(\sigma^{\cdot,l})\ast\chi_k]\,\psi_k
  +\<\sigma^{\cdot,l}\ast\chi_k,\nabla\psi_k\>.$$
It is clear that $\big|[\div(\sigma^{\cdot,l})\ast\chi_k]\,\psi_k\big|
\leq |\div(\sigma^{\cdot,l})|\ast\chi_k$. Since
  $$|\nabla\psi_k(x)|\leq \frac{\|\nabla\psi\|_\infty}{k}\,\ch_{\{k\leq|x|\leq 2k\}}
  \leq \frac{C}{1+|x|},$$
we have by \eqref{appendix-lem.1},
  $$|\<\sigma^{\cdot,l}\ast\chi_k,\nabla\psi_k\>|\leq C\frac{|\sigma^{\cdot,l}\ast\chi_k|}{1+|x|}
  \leq 2C |\bar\sigma^{\cdot,l}|\ast\chi_k.$$
Summarizing these discussions, we obtain
  \begin{equation}\label{sect-2-lem-3.1.5}
  |\div(\sigma_k^{\cdot,l})|\leq |\div(\sigma^{\cdot,l})|\ast\chi_k
  +2C|\bar\sigma^{\cdot,l}|\ast\chi_k.
  \end{equation}
Hence by Jensen's inequality,
  \begin{equation}\label{sect-2-lem-3.2}
  |\div(\sigma_k)|^2= \sum_{l=1}^m|\div(\sigma_k^{\cdot,l})|^2
  \leq 2|\div(\sigma)|^2\ast\chi_k +8C^2|\bar\sigma|^2\ast\chi_k.
  \end{equation}

Now by the definition of $\lambda$, one has $\nabla\lambda(x)=-2\alpha x/(1+|x|^2)$.
Therefore
  $$|\sigma_k^\ast\nabla\lambda|\leq 4\alpha\frac{|\sigma|\ast\chi_k}{1+|x|}
  \leq 8\alpha |\bar\sigma|\ast\chi_k,$$
where the last inequality follows from \eqref{appendix-lem.1}. As a result,
  \begin{equation}\label{sect-2-lem-3.3}
  |\sigma_k^\ast\nabla\lambda|^2\leq 64\alpha^2|\bar\sigma|^2\ast\chi_k.
  \end{equation}
Combining \eqref{sect-2-lem-3.1} with \eqref{sect-2-lem-3.2} and \eqref{sect-2-lem-3.3},
we get the estimate.

(2) Now we estimate
  $$\Lambda_2^{\sigma_k,b_k}=\div(b_k)+\L_k\lambda
  -\frac12\<\nabla\sigma_k,(\nabla\sigma_k)^\ast\>,$$
where $\L_k\lambda=\frac12\<\sigma_k\sigma_k^\ast,\textup{Hess}(\lambda)\>
+\<b_k,\nabla\lambda\>$. First we have
  $$\div(b_k)=(\div(b)\ast\chi_k)\,\psi_k+\<b\ast\chi_k,\nabla\psi_k\>,$$
and similar to the treatment of $\sigma_k^{\cdot,l}$,
  $$|\<b\ast\chi_k,\nabla\psi_k\>|\leq 2C |\bar b|\ast\chi_k .$$
Hence
  \begin{align}\label{sect-2-lem-3.4}
  [\div(b_k)]^- &\leq [\div(b)\ast\chi_k]^- +2C|\bar b|\ast\chi_k
  \leq [\div(b)]^- \ast\chi_k +2C|\bar b|\ast\chi_k.
  \end{align}

Now notice that
  $$\partial_i\partial_j\lambda(x)=-\frac{2\alpha\delta_{ij}}{1+|x|^2}
  +\frac{4\alpha x_ix_j}{(1+|x|^2)^2},$$
thus $|\partial_i\partial_j\lambda(x)|\leq C/(1+|x|)^2$ for all $x\in\R^n$.
This together with \eqref{appendix-lem.1} leads to
  \begin{align}\label{sect-2-lem-3.5}
  |\L_k\lambda|\leq C\big(|\bar\sigma|^2\ast\chi_k +|\bar b|\ast\chi_k\big).
  \end{align}
Finally, similar arguments work for estimating $\nabla\sigma_k$ and we have
  \begin{equation}\label{sect-2-lem-3.6}
  |\<\nabla\sigma_k,(\nabla\sigma_k)^\ast\>|\leq |\nabla\sigma_k|^2
  \leq C\big(|\nabla\sigma|^2\ast\chi_k +|\bar\sigma|^2\ast\chi_k \big).
  \end{equation}
Now we complete the proof by substituting the estimates \eqref{sect-2-lem-3.4}--\eqref{sect-2-lem-3.6} into
the expression of $\Lambda_2^{\sigma_k,b_k}$. \fin

\begin{lemma}[Uniform density estimate]\label{sect-2-lem-4}
For fixed $p>1$, there are two positive constants $C_{1,p},C_{2,p}>0$ and sufficiently small $T_0>0$,
such that for all $k\geq1$,
  \begin{align}\label{sect-2-lem-4.1}
  \sup_{0\leq t\leq T_0}\|\rho^k_t\|_{L^p(\P\times\mu)}
  &\leq C_{1,p}\bigg(\int_{\R^n}\exp\big[C_{2,p}T_0\big([\div(b)]^- +|\bar b|
  +|\nabla\sigma|^2+|\bar\sigma|^2\big)\big]\d \mu\bigg)^{\frac1{p(p+1)}}<+\infty.
  \end{align}
\end{lemma}

\noindent{\bf Proof.} By Lemma \ref{sect-2-lem-3}
and noticing that $|\div(\sigma)|\leq |\nabla\sigma|$, we have for any $t>0$,
  \begin{align*}
  p^3t|\Lambda_1^{\sigma_k}|^2-p^2t\Lambda_2^{\sigma_k,b_k}
  \leq Cp^3t\big[\big([\div(b)]^- +|\bar b|
  +|\nabla\sigma|^2+|\bar\sigma|^2\big)\ast\chi_k\big].
  \end{align*}
Substituting this estimate into \eqref{sect-2.2}, we see that there are two constants
$C_{1,p},C_{2,p}>0$ such that for any $T>0$ and all $k\geq1$,
  $$\sup_{0\leq t\leq T}\|\rho^k_t\|_{L^p(\P\times\mu)} \leq
  C_{1,p}\bigg(\int_{\R^n}\exp\big[C_{2,p}T\big([\div(b)]^- +|\bar b|
  +|\nabla\sigma|^2+|\bar\sigma|^2\big)\ast\chi_k\big]\d\mu\bigg)^{\frac1{p(p+1)}}.$$
To simplify the notations, we denote by $\Phi=C_{2,p}T\big([\div(b)]^- +|\bar b|
+|\nabla\sigma|^2+|\bar\sigma|^2\big)$; then
  \begin{equation}\label{sect-2-lem-4.2}
  \sup_{k\geq1}\sup_{0\leq t\leq T}\|\rho^k_t\|_{L^p(\P\times\mu)}
  \leq C_{1,p}\bigg(\int_{\R^n}\exp\big[(\Phi\ast\chi_k)(x)+\lambda(x)\big]\d x\bigg)^{\frac1{p(p+1)}}.
  \end{equation}

We want to show that there is a constant $C>0$ such that for any $k\geq1$,
  \begin{equation}\label{sect-2-lem-4.3}
  \lambda(x)\leq (\lambda\ast\chi_k)(x)+C\quad\mbox{for all } x\in\R^{n}.
  \end{equation}
Indeed, for any $u\in B(1)$, one has
  $$1+|x-u|^2\leq 1+2|x|^2+2|u|^2
  \leq 3(1+|x|^2),$$
hence
  $$\lambda(x-u)=-\alpha\log(1+|x-u|^2)
  \geq -\alpha\log 3+\lambda(x).$$
As a result, for all $k\geq 1$,
  $$(\lambda\ast\chi_k)(x)
  =\int_{\R^{n}}\lambda(x-u)\chi_k(u)\,\d u
  \geq -\alpha\log 3+\lambda(x)$$
since $\chi_k\geq0$ and $\int_{\R^{n}}\chi_k(u)\,\d u=1$. Hence \eqref{sect-2-lem-4.3}
holds with $C=\alpha\log 3$. Now by  \eqref{sect-2-lem-4.3} and  Jensen's inequality,
  \begin{align*}
  \int_{\R^n}\exp\big[(\Phi\ast\chi_k)(x)+\lambda(x)\big]\d x
  &\leq 3^\alpha \int_{\R^n}\exp\big[(\Phi+\lambda)\ast\chi_k(x)\big]\d x\cr
  &\leq 3^\alpha \int_{\R^n} (e^{\Phi+\lambda}\ast\chi_k)(x)\,\d x\cr
  &= 3^\alpha \int_{\R^n} e^{\Phi+\lambda}\,\d x
  =3^\alpha \int_{\R^n} e^\Phi\,\d\mu.
  \end{align*}
Substituting this estimate into \eqref{sect-2-lem-4.2} and by the definition
of $\Phi$, we see that if we take $T_0\leq p_0/C_{2,p}$, then the right hand side
of \eqref{sect-2-lem-4.1} is finite. \fin

\medskip

In the following we fix $p$ as the conjugate number of $q$ and denote by
$\Lambda_{p,T_0}$ the quantity on the right hand side of \eqref{sect-2-lem-4.1}.
Then we have
  \begin{equation}\label{sect-2.1}
  \sup_{k\geq1}\sup_{0\leq t\leq T_0}\|\rho^k_t\|_{L^p(\P\times\mu)}
  \leq \Lambda_{p,T_0}.
  \end{equation}
Using Lemma \ref{sect-2-lem-2} and the density estimate \eqref{sect-2.1},
we can now show that there exists a random field
$X:\Omega\times\R^n\ra C([0,T_0],\R^n)$, which is the limit
of the sequence of stochastic flows generated by \eqref{smooth-SDE}.

\begin{proposition}\label{sect-2-prop-1} Under the conditions (C1) and (C2), there
exists a random field $X:\Omega\times\R^n\ra C([0,T_0],\R^n)$ such that
  $$\lim_{k\ra\infty}\E\int_{\R^n}1\wedge\|X^k-X\|_{\infty,T_0}\,\d\mu=0.$$
\end{proposition}

\noindent{\bf Proof.} The proof is similar to that of \cite[Theorem 5.3]{FangLuoThalmaier}.
For any $k\geq1$, we denote by $G_R^{k}$ the level set of the flow $X^k_t$
on the interval $[0,T_0]$:
  $$G_R^k=\{(\omega,x)\in\Omega\times\R^n:\|X^k_\cdot(\omega,x)\|_{\infty,T_0}\leq R\}.$$
By Lemma \ref{sect-2-lem-1},
  \begin{align}\label{sect-2-prop-1.1}
  (\P\times\mu)\big[(G_R^{k}\cap G_R^l)^c\big]&\leq (\P\times\mu)\big[(G_R^k)^c\big]+(\P\times\mu)\big[(G_R^l)^c\big]
  \leq \frac{C_k+C_l}{R},
  \end{align}
in which $C_k$ depends on $T_0, \Lambda_{p,T_0}$, $\|\sigma_k\|_{L^{2q}(\mu)},
\|b_k\|_{L^{q}(\mu)}$. We have $|\sigma_k|\leq |\sigma|\ast\chi_k$.
Jensen's inequality leads to
  $$\|\sigma_k\|_{L^{2q}(\mu)}^{2q}\leq
  \int_{\R^n}\big(|\sigma|^{2q} \ast\chi_k\big)(x)\,\d\mu(x)
  =\int_{\R^n}|\sigma(y)|^{2q}\,\d y\int_{\R^n}\frac{\chi_k(x-y)}{(1+|x|^2)^\alpha}\,\d x.$$
Notice that for $|x-y|\leq 1/k$, one has $|y|\leq |x|+1/k$, hence
  $$1+|y|^2\leq 1+2|x|^2+2/k^2\leq 3(1+|x|^2)\quad \mbox{for all }k\geq1.$$
Consequently,
  \begin{equation}\label{sect-2-prop-1.1.5}
  \int_{\R^n}\frac{\chi_k(x-y)}{(1+|x|^2)^\alpha}\,\d x
  \leq 3^\alpha \int_{\R^n}\frac{\chi_k(x-y)}{(1+|y|^2)^\alpha}\,\d x
  =\frac{3^\alpha}{(1+|y|^2)^\alpha}
  \end{equation}
since $\int_{\R^n}\chi_k\,\d x=1$. As a result,
  \begin{equation}\label{sect-2-prop-1.2}
  \|\sigma_k\|_{L^{2q}(\mu)}
  \leq 3^{\alpha/2q}\bigg(\int_{\R^n}|\sigma(y)|^{2q}\,\d\mu(y)\bigg)^{1/2q}
  =3^{\alpha/2q}\|\sigma\|_{L^{2q}(\mu)}.
  \end{equation}
In the same way, we have $\|b_k\|_{L^{q}(\mu)}\leq 3^{\alpha/q}\|b\|_{L^{q}(\mu)}$.
Therefore the positive constants
$(C_k)_{k\geq1}$ are uniformly bounded from above by some $\hat C>0$.
Combining this observation with \eqref{sect-2-prop-1.1}, we obtain
  \begin{equation}\label{sect-2-prop-1.3}
  \sup_{k,l\geq 1}(\P\times\mu)\big[(G_R^k\cap G_R^l)^c\big]\leq \frac{2\hat C}{R}.
  \end{equation}

Now an application of Lemma \ref{sect-2-lem-2} to the flows $X^k_t$ and $X^l_t$ gives us
  \begin{align}\label{sect-2-prop-1.4}
  &\hskip-10pt \E\int_{G_R^k\cap G_R^l}\log\bigg(\frac{\|X^k- X^l\|^2_{\infty,T_0}}{\delta^2}+1\bigg)\d\mu\cr
  &\hskip-10pt \quad \leq C_{T_0}\Lambda_{p,T_0}\bigg\{C_{n,q}\Big[\|\nabla b_k\|_{L^q(B(3R))}
  +\|\nabla \sigma_k\|_{L^{2q}(B(3R))}+\|\nabla\sigma_k\|^2_{L^{2q}(B(3R))}\Big]\cr
  &\hskip50pt +\frac1{\delta^2}\|\sigma_k-\sigma_l\|^2_{L^{2q}(B(R))}
  +\frac1\delta \Big[\|\sigma_k-\sigma_l\|_{L^{2q}(B(R))}+
  \|b_k- b_l\|_{L^{q}(B(R))}\Big]\bigg\}.
  \end{align}
By the definition of $b_k$ and \eqref{appendix-lem.1}, we have
  $$|\nabla b_k|\leq |\nabla b|\ast\chi_k+C\frac{|b\ast\chi_k|}{1+|x|}
  \leq |\nabla b|\ast\chi_k+ 2C|\bar b|\ast \chi_k.$$
From this we can show that
  $$\|\nabla b_k\|_{L^q(B(3R))}\leq C_q\big(\|\nabla b\|_{L^q(B(3R+1))}
  +\|\bar b\|_{L^q(B(3R+1))}\big).$$
In the same way, $\|\nabla \sigma_k\|_{L^{2q}(B(3R))}\leq
C_q\big(\|\nabla \sigma\|_{L^{2q}(B(3R+1))}+\|\bar\sigma\|_{L^q(B(3R+1))}\big)$.
Notice that under the conditions (C1) and (C2), $\nabla b$ and $\bar b$ (resp.
$\nabla \sigma$ and $\bar \sigma$) are locally integrable.
Hence for any $k\geq 1$,
  $$C_{n,q}\Big[\|\nabla b_k\|_{L^q(B(3R))}+\|\nabla \sigma_k\|_{L^{2q}(B(3R))}
  +\|\nabla\sigma_k\|^2_{L^{2q}(B(3R))}\Big]  \leq C'_{n,q,R}.$$
Now we define
  $$\delta_{k,l}=\|\sigma_k-\sigma_l\|_{L^{2q}(B(R))}+\|b_k- b_l\|_{L^{q}(B(R))}$$
which tends to $0$ as $k,l\ra+\infty$. Taking $\delta=\delta_{k,l}$
in \eqref{sect-2-prop-1.4}, we obtain that for any $k,l\geq1$,
  \begin{equation}\label{sect-2-prop-1.5}
  \E\int_{G_R^k\cap G_R^l}\log\bigg(\frac{\|X^k- X^l\|^2_{\infty,T_0}}{\delta_{k,l}^2}+1\bigg)\d\mu
  \leq C_{T_0,n,q,R}<+\infty.
  \end{equation}

We have by \eqref{sect-2-prop-1.3}
  \begin{align}\label{sect-2-prop-1.6}
  &\E\int_{\R^n}\big(1\wedge\|X^k-X^l\|_{\infty,T_0}\big)\,\d\mu\cr
  &\quad \leq (\P\times\mu)\big[(G_R^k\cap G_R^l)^c\big]
  +\int_{G_R^k\cap G_R^l}\big(1\wedge\|X^k-X^l\|_{\infty,T_0}\big)\,\d(\P\times\mu)\cr
  &\quad \leq \frac{2\hat C}{R}
  +\int_{G_R^k\cap G_R^l}\big(1\wedge\|X^k-X^l\|_{\infty,T_0}\big)\,\d(\P\times\mu).
  \end{align}
Next for $\eta\in (0,1)$, set
  $$\Sigma^{k,l}_\eta=\big\{(\omega,x)\in\Omega\times\R^n:\|X^k-X^l\|_{\infty,T_0}\leq \eta\big\}.$$
Then
  \begin{align*}
  &\int_{G_R^k\cap G_R^l}\big(1\wedge\|X^k-X^l\|_{\infty,T_0}\big)\,\d(\P\times\mu)\cr
  &\quad =\bigg( \int_{(G_R^k\cap G_R^l)\cap \Sigma^{k,l}_\eta}+\int_{(G_R^k\cap G_R^l)\setminus \Sigma^{k,l}_\eta}\bigg)
  \big(1\wedge\|X^k-X^l\|_{\infty,T_0}\big)\,\d(\P\times\mu)\cr
  &\quad \leq \eta\,\mu(\R^n)+\frac{1}{\log\Big(1+\frac{\eta^2}{\delta_{k,l}^2}\Big)}
  \int_{G_R^k\cap G_R^l}\log\bigg(1+\frac{\|X^k-X^l\|_{\infty,T_0}^2}{\delta_{k,l}^2}\bigg)\d(\P\times\mu)\cr
  &\quad \leq \eta\,\mu(\R^n)+\frac{C_{T_0,n,q,R}}{\log\Big(1+\frac{\eta^2}{\delta_{k,l}^2}\Big)},
  \end{align*}
where the last inequality follows from \eqref{sect-2-prop-1.5}. Substituting
this estimate into \eqref{sect-2-prop-1.6}, we get
  \begin{equation*}
  \E\int_{\R^n}\big(1\wedge\|X^k-X^l\|_{\infty,T_0}\big)\,\d\mu
  \leq \frac{2\hat C}{R} +\eta\,\mu(\R^n)
  +\frac{C_{T_0,n,q,R}}{\log\Big(1+\frac{\eta^2}{\delta_{k,l}^2}\Big)}.
  \end{equation*}
First letting $k,l\ra+\infty$, and then $R\ra+\infty$, $\eta\ra0$, we obtain
  $$\lim_{k,l\ra+\infty}\E\int_{\R^n}\big(1\wedge\|X^k-X^l\|_{\infty,T_0}\big)\,\d\mu=0.$$
Hence there exists a random field $X:\Omega\times\R^n\ra C([0,T_0],\R^n)$ such that
  \begin{equation*}
  \lim_{k\ra+\infty}\E\int_{\R^n}\big(1\wedge\|X^k-X\|_{\infty,T_0}\big)\,\d\mu=0.
  \end{equation*}

\begin{proposition}\label{sect-2-prop-2} For all $t\in[0,T_0]$, there exists
$\rho_t:\Omega\times\R^n\ra\R_+$ such that
$(X_t)_\#\mu=\rho_t \mu$. Moreover, $\sup_{0\leq t\leq
T_0}\|\rho_t\|_{L^p(\P\times\mu)}\leq \Lambda_{p,T_0}$.
\end{proposition}

\noindent{\bf Proof.} We follow the arguments of \cite[Theorem 3.4]{FangLuoThalmaier}.
By Proposition \ref{sect-2-prop-1}, it is easy to show that
for any $\psi\in C^\infty_c(\R^n)$
  $$\lim_{k\ra\infty}\E\int_{\R^n}|\psi(X^k_t(x))-\psi(X_t(x))|\,\d\mu=0.$$
Now we fix any $\xi\in L^\infty(\Omega)$ and $\psi\in
C^\infty_c(\R^n)$; then
  $$\E\int_{\R^n}|\xi(\omega)|\cdot|\psi(X^k_t(x))-\psi(X_t(x))|\,\d\mu
  \leq\|\xi\|_\infty\E\int_{\R^n}|\psi(X^k_t(x))-\psi(X_t(x))|\,\d\mu \ra0$$
as $k$ goes to $\infty$. Thus we have
  $$\lim_{k\ra\infty}\E\int_{\R^n}\xi(\omega)\psi(X^k_t(x))\,\d\mu
  =\E\int_{\R^n}\xi(\omega)\psi(X_t(x))\,\d\mu.$$
On the other hand, since $(X^k_t)_\#\mu=\rho^k_t\mu$ and the
family $\{\rho^k_t:k\geq1\}$ is bounded in $L^p(\Omega\times\R^n)$ for
all $t\leq T_0$, thus up to a subsequence, $\rho^k_t$ converges weakly
to some $\rho_t\in L^p(\Omega\times\R^n)$. By the property of
weak convergence, we have
  $$\|\rho_t\|_{L^p(\P\times\mu)}
  \leq \Lambda_{p,T_0},\quad\mbox{for all }t\leq T_0.$$
Therefore
  $$\lim_{k\ra\infty}\E\int_{\R^n}\xi(\omega)\psi(X^k_t(x))\,\d\mu
  =\lim_{k\ra\infty}\E\int_{\R^n}\xi(\omega)\psi(y)\rho^k_t(y)\,\d\mu
  =\E\int_{\R^n}\xi(\omega)\psi(y)\rho_t(y)\,\d\mu.$$
Combining the above two equalities, we obtain for all $t\leq T_0$,
  $$\E\int_{\R^n}\xi(\omega)\psi(X_t(x))\,\d\mu
  =\E\int_{\R^n}\xi(\omega)\psi(y)\rho_t(y)\,\d\mu.$$
By the arbitrariness of $\xi\in L^\infty(\Omega)$, there is a full
subset $\Omega_\psi$ such that for all $\omega\in\Omega_\psi$, it
holds
  $$\int_{\R^n}\psi(X_t(x))\,\d\mu=\int_{\R^n}\psi(y)\rho_t(y)\,\d\mu,
  \quad\mbox{for all }\psi\in C_c^\infty(\R^n).$$
Now by the separability of $C_c^\infty(\R^n)$, we can find another
full subset $\Omega_t$, such that for every $\omega\in \Omega_t$,
the above equality holds for all $\psi\in C_c^\infty(\R^n)$. From
this we conclude that $(X_t)_\#\mu=\rho_t\mu$.  \fin

\medskip

To show that $(X_t)_{0\leq t\leq T_0}$ solves the It\^o SDE \eqref{Ito-SDE}, we need
the following preparations.

\begin{lemma}\label{sect-2-lem-5} We have
  $$\lim_{k\ra\infty}\|\sigma_k-\sigma\|_{L^{2q}(\mu)}=0
  \quad\mbox{and}\quad \lim_{k\ra\infty}\|b_k-b\|_{L^{2q}(\mu)}=0.$$
\end{lemma}

\noindent{\bf Proof.} By the triangular inequality,
  \begin{equation}\label{sect-2-lem-5.1}
  \|\sigma_k-\sigma\|_{L^{2q}(\mu)}\leq \|(\sigma\ast\chi_k)(\psi_k-1)\|_{L^{2q}(\mu)}
  +\|\sigma\ast\chi_k-\sigma\|_{L^{2q}(\mu)}.
  \end{equation}
We deduce from Jensen's inequality that
  $$\int_{\R^n}|(\sigma\ast\chi_k)(\psi_k-1)|^{2q}\,\d\mu
  \leq \int_{\R^n}(1-\psi_k)|\sigma\ast\chi_k|^{2q}\,\d\mu
  \leq \int_{\R^n}\ch_{\{|x|\geq k\}}|\sigma|^{2q}\ast\chi_k\,\d\mu.$$
Fubini's theorem leads to
  \begin{align*}
  \int_{\R^n}\ch_{\{|x|\geq k\}}|\sigma|^{2q}\ast\chi_k\,\d\mu
  &= \int_{\R^n}\!\int_{\R^n}|\sigma(y)|^{2q}
  \frac{\chi_k(x-y)\ch_{\{|x|\geq k\}}}{(1+|x|^2)^\alpha}\,\d y\d x\cr
  &\leq \int_{\R^n}|\sigma(y)|^{2q}\ch_{\{|y|\geq k-1\}}\,\d y
  \int_{\R^n}\frac{\chi_k(x-y)}{(1+|x|^2)^\alpha}\,\d x.
  \end{align*}
Thus by \eqref{sect-2-prop-1.1.5}, we obtain
  $$\int_{\R^n}|(\sigma\ast\chi_k)(\psi_k-1)|^{2q}\,\d\mu
  \leq 3^\alpha\int_{\R^n}|\sigma(y)|^{2q}\ch_{\{|y|\geq k-1\}}\,\d\mu(y).$$
Notice that $\sigma\in L^{2q}(\mu)$ (see Remark \ref{sect-2-rem}), we deduce that
  \begin{equation}\label{sect-2-lem-5.2}
  \lim_{k\ra\infty}\|(\sigma\ast\chi_k)(\psi_k-1)\|_{L^{2q}(\mu)}=0.
  \end{equation}

Next for any $R>0$, we have
  \begin{equation*}
  \int_{\R^n}|\sigma\ast\chi_k-\sigma|^{2q}\,\d\mu
  = \bigg(\int_{\{|x|\leq R\}}+\int_{\{|x|> R\}}\bigg)|\sigma\ast\chi_k-\sigma|^{2q}\,\d\mu.
  \end{equation*}
By the above discussions, it is clear that
  \begin{align*}
  \int_{\{|x|> R\}}|\sigma\ast\chi_k-\sigma|^{2q}\,\d\mu
  &\leq C_q\int_{\{|x|> R\}}|\sigma\ast\chi_k|^{2q}\,\d\mu+C_q\int_{\{|x|> R\}}|\sigma|^{2q}\,\d\mu\cr
  &\leq 3^\alpha C_q\int_{\{|x|> R-1\}}|\sigma|^{2q}\,\d\mu+C_q\int_{\{|x|> R\}}|\sigma|^{2q}\,\d\mu.
  \end{align*}
Thus for any $k\geq1$,
  \begin{equation}\label{sect-2-lem-5.3}
  \int_{\R^n}|\sigma\ast\chi_k-\sigma|^{2q}\,\d\mu\leq
  \int_{\{|x|\leq R\}}|\sigma\ast\chi_k-\sigma|^{2q}\,\d\mu
  +(3^\alpha+1) C_q\int_{\{|x|> R-1\}}|\sigma|^{2q}\,\d\mu.
  \end{equation}
It is obvious that
  $$\lim_{k\ra\infty}\int_{\{|x|\leq R\}}|\sigma\ast\chi_k-\sigma|^{2q}\,\d\mu
  \leq \lim_{k\ra\infty}\int_{\{|x|\leq R\}}|\sigma\ast\chi_k-\sigma|^{2q}\,\d x=0$$
for any fixed $R>0$. Hence first letting $k\ra\infty$ and then $R\ra\infty$
in \eqref{sect-2-lem-5.3}, we obtain
  $$\lim_{k\ra\infty}\int_{\R^n}|\sigma\ast\chi_k-\sigma|^{2q}\,\d\mu=0.$$
Combining this with \eqref{sect-2-lem-5.1} \eqref{sect-2-lem-5.2}, we
obtain the first result. The second one can be proved in the same way, hence we omit it. \fin

\begin{corollary}\label{sect-2-cor-2} We have
  $$\lim_{k\ra\infty}\E\int_{\R^n}\bigg(\sup_{0\leq t\leq T_0}
  \bigg|\int_0^t\big[\sigma_k(X^k_s)-\sigma(X_s)\big]\d B_s\bigg|\bigg)\d\mu=0$$
and
  $$\lim_{k\ra\infty}\E\int_{\R^n}\bigg(\sup_{0\leq t\leq T_0}
  \bigg|\int_0^t\big[b_k(X^k_s)-b(X_s)\big]\d s\bigg|\bigg)\d\mu=0.  $$
\end{corollary}

\noindent{\bf Proof.} Having Propositions \ref{sect-2-prop-1}, \ref{sect-2-prop-2}
and Lemma \ref{sect-2-lem-5} in mind, the proof is similar to that of
\cite[Proposition 4.1]{FangLuoThalmaier}. We omit it here. \fin

\medskip

For any $k\geq1$, we rewrite the equation \eqref{smooth-SDE} in the integral form:
  \begin{equation}\label{smooth-SDE-integral}
  X^k_t(x)=x+\int_0^t\sigma_k(X^k_s)\,\d B_s
  +\int_0^t b_k(X^k_s)\,\d s.
  \end{equation}
When $k\ra +\infty$, by Proposition \ref{sect-2-prop-1} and Corollary
\ref{sect-2-cor-2}, the two sides of \eqref{smooth-SDE-integral} converge respectively
to $X$ and
  \begin{equation*}
  x+\int_0^{\textstyle\cdot} \sigma(X_s)\,\d B_s
  +\int_0^{\textstyle\cdot} b(X_s)\,\d s.
  \end{equation*}
Therefore, for almost all $x\in\R^d$, the following equality holds
$\P$-almost surely:
  $$X_t(x)=x+\int_0^t \sigma(X_s)\,\d B_s
  +\int_0^t b(X_s)\,\d s,\quad\mbox{for all }t\in[0,T_0].$$
That is to say, $X_t$ solves SDE \eqref{Ito-SDE} over the time interval $[0,T_0]$.
Similar to \cite[Proposition 5.6]{FangLuoThalmaier}, we can prove the uniqueness
of the solution flow on $[0,T_0]$.

Now we extend the solution to any time interval $[0,T]$. Let
$\theta_{T_0}B$ be the time-shift of the Brownian motion $B$ by $T_0$ and
denote by $X_t^{T_0}$ the corresponding solution to the SDE \eqref{Ito-SDE} driven by
$\theta_{T_0}B$. By the above discussions,
$\{X_t^{T_0}(\theta_{T_0}B,x)\colon\, 0\leq t\leq T_0\}$ is the unique
solution to the following SDE over $[0,T_0]$:
  \begin{equation*}
  X_t^{T_0}(x)=x+ \int_0^t \sigma \big(X_s^{T_0}(x)\big)\,\d(\theta_{T_0}B)_s
  +\int_0^t b\big(X_s^{T_0}(x)\big)\,\d s.
  \end{equation*}
For $t\in[0,T_0]$, define $X_{t+T_0}(\omega,x)=X^{T_0}_t(\theta_{T_0}B,X_{T_0}(\omega,x))$.
Note that $X_t$ is well defined on the interval $[0,2T_0]$ up to a
$(\P\times\mu)$-negligible subset of $\Omega\times\R^n$. Replacing $x$ by $X_{T_0}(x)$
in the above equation, we obtain
  \begin{equation*}
  X_{t+T_0}(x)=x+ \int_0^{t+T_0} \sigma(X_s(x))\,\d B_s
  +\int_0^{t+T_0} b(X_s(x))\,\d s.
  \end{equation*}
Therefore $X_t$ defined as above is a solution to SDE \eqref{Ito-SDE} on the
interval $[0,2T_0]$. Continuing in this way, we obtain the solution
of SDE \eqref{Ito-SDE} on the interval $[0,T]$.

\begin{proposition}\label{sect-2-prop-3}
The family $\{X_t\colon t\in [0,T]\}$ constructed as above is the unique solution
to SDE \eqref{Ito-SDE}.
\end{proposition}

\noindent{\bf Proof.} Let $Y_t$, $t\in[0,T]$ be another solution. First by the above discussions,
we have $(\P\times \mu)$-almost surely,
$Y_t=X_t$ for all $t\in [0,T_0]$. In particular, $Y_{T_0}=X_{T_0}$.
Next by the flow property, $Y_{t+T_0}$ satisfies the following
equation:
  \begin{equation*}
  Y_{t+T_0}(x)=Y_{T_0}(x)+\int_0^t \sigma\big(Y_{s+T_0}(x)\big)\d(\theta_{T_0}B)_s
  +\int_0^t b\big(Y_{s+T_0}(x)\big)\d s,
  \end{equation*}
that is, $Y_{t+T_0}$ is a solution with initial value $Y_{T_0}$. But
by the above discussion, $X_{t+T_0}$ is also a solution with the
same initial value $X_{T_0}=Y_{T_0}$. Therefore, we have $(\P\times\mu)$-almost surely,
$X_{t+T_0}=Y_{t+T_0}$ for all $t\leq T_0$. Hence we have proved that
$X|_{[0,2T_0]}=Y|_{[0,2T_0]}$. Repeating this procedure, we obtain
the uniqueness over $[0,T]$. \fin

\medskip

Now we want to show that the reference measure $\mu$ is absolutely continuous
under the stochastic flow $\{X_t:t\leq T\}$ constructed above. To this end
we have to prove an $L^1\log L^1$-type estimate for the density functions $\rho^k_t$,
and extend the convergence result in Proposition \ref{sect-2-prop-1}
to general time interval $[0,T]$.

\begin{proposition}\label{sect-2-prop-5}
For each $t\in [0,T]$, there exists $\rho_t:\Omega\times\R^n\ra\R_+$ such that
$\rho_t\in L^1\log L^1$.
\end{proposition}

\noindent{\bf Proof.} Let $T>0$ be given. Similar to \cite[Theorem 3.3]{FangLuoThalmaier},
we can prove
  \begin{equation}\label{sect-2-prop-5.1}
  \sup_{k\geq1}\sup_{0\leq t\leq T}\E\int_{\R^n}\rho^k_t|\log\rho^k_t|\,\d\mu<+\infty.
  \end{equation}
Here we give a sketch of the proof. We have
  $$\int_{\R^n}\rho^k_t|\log\rho^k_t|\,\d\mu=\int_{\R^n}\big|\log\rho^k_t(X^k_t)\big|\,\d\mu.$$
By \eqref{appendix.1.5}, $\rho^k_t(X^k_t)=1/\tilde\rho^k_t$, where $\tilde\rho^k_t$ is the
density of $\big((X^k_t)^{-1}\big)_\#\mu$ with respect to $\mu$ which admits the expression
\eqref{appendix.2}. Using the flow property and \eqref{sect-2.1}, the It\^o calculus
leads to \eqref{sect-2-prop-5.1}. Using the estimate \eqref{sect-2-prop-5.1}, for each
$t\in[0,T]$, there is a subsequence $k_i$ such that $\rho^{k_i}_t$ converge weakly in
$L^1(\Omega\times\R^n)$ to $\rho_t$. Following the argument on the page 1144 of \cite{FangLuoThalmaier},
we conclude that $\rho_t\in L^1\log L^1$.  \fin

\medskip

\begin{proposition}\label{sect-2-prop-4}
For any $t\in[0,T]$, $\rho_t$ is the density of $(X_t)_\#\mu$ with respect to
$\mu$; moreover we have
  $$\lim_{k\ra\infty}\E\int_{\R^n}1\wedge\|X^k-X\|_{\infty,T}\,\d\mu=0.$$
\end{proposition}

\noindent{\bf Proof.} We shall prove this result by induction. By Propositions \ref{sect-2-prop-1} and
\ref{sect-2-prop-2}, we see that the assertions are true on the interval $[0,T_0]$.
Now suppose we have proved the assertions on the time interval $[0,lT_0]$ where $l\geq1$.

For $t\leq T_0$, using the flow property of $X^k_t$ and $X_t$, we have
  \begin{align*}
  \big|X^k_{t+lT_0}-X_{t+lT_0}\big|&=\big|X^{k,lT_0}_t(X^k_{lT_0})-X^{lT_0}_t(X_{lT_0})\big|\cr
  &\leq \big|X^{k,lT_0}_t(X^k_{lT_0})-X^{lT_0}_t(X^k_{lT_0})\big|
  +\big|X^{lT_0}_t(X^k_{lT_0})-X^{lT_0}_t(X_{lT_0})\big|,
  \end{align*}
where $X^{k,lT_0}_t$ \big(resp. $X^{lT_0}_t$\big) is the solution of \eqref{smooth-SDE}
(resp. \eqref{Ito-SDE}) driven by the shifted Brownian motion $\theta_{lT_0}B$. Therefore
  \begin{align}\label{sect-2-prop-4.1}
  \E\int_{\R^n}1\wedge \big\|X^k_{\cdot+lT_0}-X_{\cdot+lT_0}\big\|_{\infty,T_0}\,\d\mu
  &\leq \E\int_{\R^n}1\wedge\big\|X^{k,lT_0}_\cdot(X^k_{lT_0})-X^{lT_0}_\cdot(X^k_{lT_0})\big\|_{\infty,T_0}\,\d\mu\cr
  &\quad + \E\int_{\R^n}1\wedge\big\|X^{lT_0}_\cdot(X^k_{lT_0})-X^{lT_0}_\cdot(X_{lT_0})\big\|_{\infty,T_0}\,\d\mu\cr
  &=: I^k_1+I^k_2.
  \end{align}
By the induction hypothesis, for any $R>1$, we have
  \begin{align}\label{sect-2-prop-4.2}
  I^k_1&=\E\int_{\R^n}\Big(1\wedge\big\|X^{k,lT_0}_\cdot(y)-X^{lT_0}_\cdot(y)\big\|_{\infty,T_0}\Big)
  \rho^k_{lT_0}(y)\,\d\mu(y)\cr
  &=\bigg(\int_{\{\rho^k_{lT_0}\leq R\}}+\int_{\{\rho^k_{lT_0}> R\}}\bigg)
  \Big(1\wedge\big\|X^{k,lT_0}_\cdot-X^{lT_0}_\cdot\big\|_{\infty,T_0}\Big)
  \rho^k_{lT_0}\,\d(\P\times\mu)\cr
  &\leq R\int_{\Omega\times\R^n}1\wedge\big\|X^{k,lT_0}_\cdot-X^{lT_0}_\cdot\big\|_{\infty,T_0}\,\d(\P\times\mu)
  +\int_{\{\rho^k_{lT_0}> R\}}\rho^k_{lT_0}\,\d(\P\times\mu).
  \end{align}
\eqref{sect-2-prop-5.1} tells us that
  $$\int_{\{\rho^k_{lT_0}> R\}}\rho^k_{lT_0}\,\d(\P\times\mu)
  \leq \frac1{\log R}\int_{\Omega\times\R^n}\rho^k_{lT_0}\big|\log\rho^k_{lT_0}\big|\,\d(\P\times\mu)
  \leq \frac C{\log R},$$
where $C>0$ is independent of $k\geq1$. Therefore by \eqref{sect-2-prop-4.2} and
Proposition \ref{sect-2-prop-1}, $\limsup_{k\ra\infty} I^k_1\leq \frac C{\log R}$.
Letting $R\ra+\infty$, we get
  \begin{align}\label{sect-2-prop-4.2.5}
  \lim_{k\ra\infty} I^k_1=0.
  \end{align}

Now we deal with the term $I^k_2$. Fix an arbitrary $\eta>0$. It is easy to see
that the results in Proposition \ref{sect-2-prop-2} still hold for $\{X^{lT_0}_t:0\leq t\leq T_0\}$.
Therefore we can apply \eqref{sect-2-lem-1-4} to get
  $$\E\int_{\R^n} \big\|X^{lT_0}_\cdot(x)\big\|_{\infty,T_0}\,\d\mu(x)<\infty;$$
that is, $X^{lT_0}_\cdot\in L^1\big(\Omega\times\R^n,C([0,T_0],\R^n)
\big)$. Hence for any $\ee>0$, there exists $\Phi=\sum_{j=1}^h\xi_j\varphi_j$
with $\xi_j\in L^\infty(\Omega)$ and $\varphi_j\in C_c(\R^n,C([0,T_0],
\R^n))$ such that
  \begin{equation}\label{sect-2-prop-4.3}
  \E\int_{\R^n}\big\|X^{lT_0}_\cdot(x)-\Phi(\cdot,x)\big\|_{\infty,T_0}\,\d\mu(x)<\ee.
  \end{equation}
By the triangular inequality,
  \begin{align*}
  I^k_2&\leq \E\int_{\R^n}1\wedge\big\|X^{lT_0}_\cdot(X^k_{lT_0})
  -\Phi(\cdot,X^k_{lT_0})\big\|_{\infty,T_0}\,\d\mu
  +\E\int_{\R^n}1\wedge\big\|\Phi(\cdot,X^k_{lT_0})-\Phi(\cdot,X_{lT_0})\big\|_{\infty,T_0}\,\d\mu\cr
  &\quad +\E\int_{\R^n}1\wedge\big\|\Phi(\cdot,X_{lT_0})-X^{lT_0}_\cdot(X_{lT_0})\big\|_{\infty,T_0}\,\d\mu\cr
  &=:I^k_{2,1}+I^k_{2,2}+I^k_{2,3}.
  \end{align*}
Analogous treatment of $I^k_1$ leads to
  \begin{align*}
  I^k_{2,1}&=\E\int_{\R^n}\Big(1\wedge\big\|X^{lT_0}_\cdot(y)
  -\Phi(\cdot,y)\big\|_{\infty,T_0}\Big) \rho^k_{lT_0}(y)\,\d\mu(y)\cr
  &= \bigg(\int_{\{\rho^k_{lT_0}\leq R\}}+\int_{\{\rho^k_{lT_0}> R\}}\bigg)
  \Big(1\wedge\big\|X^{lT_0}_\cdot(y)-\Phi(\cdot,y)\big\|_{\infty,T_0}\Big)
  \rho^k_{lT_0}(y)\,\d(\P\times\mu)\cr
  &\leq R \int_{\Omega\times\R^n}1\wedge\big\|X^{lT_0}_\cdot(y)-
  \Phi(\cdot,y)\big\|_{\infty,T_0}\,\d(\P\times\mu)
  +\frac1{\log R}\int_{\Omega\times\R^n}\rho^k_{lT_0}\big|\log\rho^k_{lT_0}\big|\,\d(\P\times\mu).
  \end{align*}
By \eqref{sect-2-prop-4.3} and \eqref{sect-2-prop-5.1}, we have
  $$I^k_{2,1}\leq R\ee+\frac C{\log R}.$$
Taking $R=e^{2C/\eta}$ and $\ee\leq \eta/2R$, we see that $I^k_{2,1}\leq \eta$
for all $k\geq1$. By the induction hypotheses,
we have $(X_{lT_0})_\#\mu=\rho_{lT_0}\,\mu$ and $X^k_{lT_0}$
converge to $X_{lT_0}$ in the measure $\P\times\mu$. Noticing that
$\rho_{lT_0}\in L^1\log L^1$, we can estimate the term $I^k_{2,3}$ in the
same way as $I^k_{2,1}$ and obtain $I^k_{2,3}\leq \eta$. Moreover,
by the dominated convergence theorem, $I^k_{2,2}\ra0$ as $k\ra\infty$.
To sum up,
  $$\limsup_{k\ra\infty}I^k_2\leq 2\eta.$$
Since $\eta>0$ is arbitrary,
this together with \eqref{sect-2-prop-4.1} and \eqref{sect-2-prop-4.2.5} leads to
  $$\lim_{k\ra\infty}\E\int_{\R^n}1\wedge
  \big\|X^k_{\cdot+lT_0}-X_{\cdot+lT_0}\big\|_{\infty,T_0}\,\d\mu=0.$$
Noticing that
  $$1\wedge\|X^k-X\|_{\infty,(l+1)T_0}\leq 1\wedge\|X^k-X\|_{\infty,lT_0}
  +1\wedge\big\|X^k_{\cdot+lT_0}-X_{\cdot+lT_0}\big\|_{\infty,T_0},$$
we conclude
  $$\lim_{k\ra\infty}\E\int_{\R^n}1\wedge\|X^k-X\|_{\infty,(l+1)T_0}\,\d\mu=0.$$
Now for any $t\in[lT_0,(l+1)T_0]$, since $\rho_t$ is the weak limit in
$L^1(\Omega\times\R^n)$ of some subsequence of $\rho^k_t$, we can repeat the
proof of Proposition \ref{sect-2-prop-2} to show that $\rho_t$ is the density
of $(X_t)_\#\mu$ with respect to $\mu$. Therefore we have proved the assertions
on the time interval $[0,(l+1)T_0]$. By the induction method, we finally get the desired result. \fin

\section{An intermediate result}

In this section we prove a technical result which serves as a bridge between
Theorem \ref{Extension} and the main result in Section 4. First we introduce some
notations. The functions $\sigma_i$ and $b_i$ ($i=1,2$) are the same as in the introduction.
Again we fix some $q>1$ and choose $\alpha_1>q+n_1/2,\, \alpha>\alpha_1+ n_2/2$. Let
  $$\d\mu(x)=(1+|x|^2)^{-\alpha}\,\d x\quad \mbox{and}\quad
  \d\mu_1(x_1)=(1+|x_1|^2)^{-\alpha_1}\,\d x_1.$$
Then $\mu$ (resp. $\mu_1$) is a finite measure on $\R^n$ (resp. $\R^{n_1}$).
To simplify the notations we write $\bar \sigma_1=\frac{\sigma_1}{1+|x_1|}$ and
$\bar \sigma_2=\frac{\sigma_2}{1+|x|}$. $\bar b_i$ is defined similarly to
$\bar\sigma_i\ (i=1,2)$.
Our assumptions in this section are:
  \begin{itemize}
  \item[\rm(H1)] $\sigma_1\in W^{1,2q}_{x_1,loc},b_1\in W^{1,q}_{x_1,loc}$;
  \item[\rm(H2)] $\int_{\R^{n_1}}\exp\big[p_0\big([\div_{x_1}(b_1)]^- +|\bar b_1|+|\bar \sigma_1|^2
  +|\nabla_{x_1}\sigma_1|^2\big)\big]\d\mu_1<+\infty$ for some $p_0>0$;
  \item[\rm(H3)] $\sigma_2\in W^{1,2q}_{x_1,x_2,loc},b_2\in W^{1,q}_{x_1,x_2,loc}$;
  \item[\rm(H4)] $\int_{\R^{n}}\exp\big[p_0\big([\div_{x_2}(b_2)]^- +|\bar b_2|+|\bar\sigma_2|^2
  +|\nabla_{x_2}\sigma_2|^2\big)\big]\d\mu<+\infty$ for some $p_0>0$.
  \end{itemize}

Under the conditions (H1) and (H2), we conclude from Theorem \ref{Extension} that there
exists a unique stochastic flow $X_{1,t}$ on $\R^{n_1}$ associated to the It\^o
SDE \eqref{Ito-SDE} with coefficients $\sigma_1$ and $b_1$, such that the reference measure $\mu_1$
is absolutely continuous under the action of the flow $X_{1,t}$. In the next result
we show that under the additional assumptions (H3)--(H4), the following SDE
  \begin{equation}\label{sect-3.1}
  \begin{cases}
  \d X_{1,t}=\sigma_1(X_{1,t})\,\d B_t+ b_1(X_{1,t})\,\d t, & X_{1,0}=x_1,\\
  \d X_{2,t}=\sigma_2(X_{1,t},X_{2,t})\,\d B_t+ b_2(X_{1,t},X_{2,t})\,\d t, & X_{2,0}=x_2
  \end{cases}
  \end{equation}
generates a unique flow $X_t=(X_{1,t},X_{2,t})$ on the whole space $\R^n$,
which leaves the measure $\mu$ absolutely continuous. Notice that the hypotheses (H1) and (H3)
imply $\sigma=(\sigma_1,\sigma_2)\in W^{1,2q}_{x_1,x_2,loc}$ and $b=(b_1,b_2)
\in W^{1,q}_{x_1,x_2,loc}$, therefore the following theorem can essentially be seen as
a special case of Theorem \ref{Extension} (see also \cite[Theorem 2.4]{Zhang12}
and \cite[Theorem 1.3]{FangLuoThalmaier}). The main difference between the two results
is that we no longer require the exponential integrability of all the
partial derivatives of $\sigma_2$; the reason for this will become clear in
view of \eqref{key-observation}.

\begin{theorem}\label{auxiliary-theorem}
Under the assumptions (H1)--(H4), the It\^o SDE \eqref{sect-3.1} generates a unique stochastic
flow $X_t$ of measurable maps on $\R^n$. Moreover, the Radon--Nikodym density $\rho_t$ of the
flow with respect to the measure $\mu$ satisfies $\rho_t\in L^1\log L^1$.
\end{theorem}

We shall not give a complete proof to the above result, but only mention some
arguments that are different from those in Section 2.
To prove Theorem \ref{auxiliary-theorem}, we need the estimates of the level sets
$G_{i,R}=\{(\omega,x):\|X_{i,\cdot}\|_{\infty,T}\leq R\}$ for the
process $X_{i,t}\ (i=1,2)$ which are similar to Lemma \ref{sect-2-lem-1}.
Notice that we do not distinguish the norms of $C([0,T],\R^n)$ and $C([0,T],\R^{n_i})
\ (i=1,2)$.

\begin{lemma}\label{sect-3-lem-1}
Let $X_t=(X_{1,t},X_{2,t})$ be a generalized stochastic flow associated to It\^o SDE
\eqref{sect-3.1}. Denote by $\rho_t$ (resp. $\rho_{1,t}$) the Radon--Nikodym density
of $X_t$ (resp. $X_{1,t}$) with respect to $\mu$ (resp. $\mu_1$). Suppose that
  $$\Lambda_{p,T}:=\sup_{0\leq t\leq T}\|\rho_t\|_{L^p(\P\times\mu)}
  \vee\|\rho_{1,t}\|_{L^p(\P\times\mu_1)} <+\infty.$$
Then under the conditions (H2) and (H4), we have
  $$(\P\times\mu_1)(G_{1,R}^c)\leq \frac{C_1}{R}\quad
  \mbox{and}\quad (\P\times\mu)(G_{2,R}^c)\leq \frac{C_2}{R},$$
where $C_1$ (resp. $C_2$) depends on $T,\Lambda_{p,T}$, $\|\sigma_1\|_{L^{2q}(\mu_1)}$ and
$\|b_1\|_{L^{q}(\mu_1)}$ (resp. $\|\sigma_2\|_{L^{2q}(\mu)}$ and
$\|b_2\|_{L^{q}(\mu)}$).

Consequently, $(\P\times\mu)(G_R^c)\leq \tilde C/R$, where $G_R$ is the level set of
$X_t=(X_{1,t},X_{2,t})$.
\end{lemma}

\noindent{\bf Proof.} First we remark that under the condition (H2) (resp. (H4)),
the coefficients $\sigma_1$ and $b_1$ (resp. $\sigma_2$ and $b_2$) belong to the space
$L^{2q}(\mu_1)$ (resp. $L^{2q}(\mu)$), see Remark \ref{sect-2-rem}(ii) for the proof.
The first estimate has been proved in Lemma \ref{sect-2-lem-1}. Here
we give a proof of the second one. We have
  \begin{equation}\label{sect-3-lem-1.1}
  \|X_{2,\cdot}\|_{\infty,T}\leq |x_2|
  +\sup_{0\leq t\leq T}\bigg|\int_0^t \sigma_2(X_s)\,\d B_s\bigg|
  +\sup_{0\leq t\leq T}\bigg|\int_0^t b_2(X_s)\,\d s\bigg|.
  \end{equation}
By Burkholder's inequality,
  $$\E\sup_{0\leq t\leq T}\bigg|\int_0^t \sigma_2(X_s)\,\d B_s\bigg|
  \leq 2\bigg[\E\int_0^T |\sigma_2(X_s)|^2\,\d s\bigg]^{\frac12}.$$
Cauchy's inequality leads to
  \begin{align*}
  \int_{\R^n}\E\sup_{0\leq t\leq T}\bigg|\int_0^t \sigma_2(X_s)\,\d B_s\bigg|\d\mu
  \leq 2\mu(\R^n)^{\frac12}\bigg[\int_0^T \E\int_{\R^n}|\sigma_2(X_s)|^2\,\d\mu\d s\bigg]^{\frac12}.
  \end{align*}
We have by H\"older's inequality that
  $$\E\int_{\R^n}|\sigma_2(X_s)|^2\,\d\mu=\E\int_{\R^n}|\sigma_2|^2\rho_s\,\d\mu
  \leq \|\sigma_2\|_{L^{2q}(\mu)}^2\|\rho_s\|_{L^p(\P\times\mu)}\leq
  \Lambda_{p,T}\|\sigma_2\|_{L^{2q}(\mu)}^2.$$
Therefore
  \begin{equation}\label{sect-3-lem-1.2}
  \int_{\R^n}\E\sup_{0\leq t\leq T}\bigg|\int_0^t \sigma_2(X_s)\,\d B_s\bigg|\d\mu
  \leq 2\big(\mu(\R^n)\,T\Lambda_{p,T}\big)^{\frac12}\|\sigma_2\|_{L^{2q}(\mu)}.
  \end{equation}
In the same way, we have
  $$\E\int_{\R^n}\sup_{0\leq t\leq T}\bigg|\int_0^t b_2(X_s)\,\d s\bigg|\d\mu
  \leq T\Lambda_{p,T}\|b_2\|_{L^{q}(\mu)}.$$
Combining this with \eqref{sect-3-lem-1.2} and integrating both sides of \eqref{sect-3-lem-1.1},
we get
  $$\E\int_{\R^n}\|X_{2,\cdot}\|_{\infty,T}\,\d\mu\leq \tilde C+
  2\big(\mu(\R^n)\,T\Lambda_{p,T}\big)^{\frac12}\|\sigma_2\|_{L^{2q}(\mu)}
  +T\Lambda_{p,T}\|b_2\|_{L^{q}(\mu)},$$
where $\tilde C=\int_{\R^n}|x_2|\,\d\mu(x)<+\infty$. Now the second estimate
of level sets follows from Chebyshev's inequality.

The last assertion is obvious from the observation below:
  \begin{align*}
  (\P\times\mu)(G_R^c)&\leq (\P\times\mu)(G_{1,R/2}^c)+(\P\times\mu)(G_{2,R/2}^c)\cr
  &\leq \mu_2(\R^{n_2})\, (\P\times\mu_1)(G_{1,R/2}^c)+(\P\times\mu)(G_{2,R/2}^c),
  \end{align*}
where $\d\mu_2(x_2)=(1+|x_2|^2)^{\alpha_1-\alpha}\,\d x_2$ is a finite measure on
$\R^{n_2}$. \fin

\medskip

Next we shall present a stability estimate of the form  Lemma \ref{sect-2-lem-2}.
Suppose we are given a matrix-valued function
$\tilde\sigma:\R^n\ra\R^m\otimes\R^n$ which has the same structure with $\sigma$,
that is $\tilde\sigma=(\tilde\sigma_1,\tilde\sigma_2)$ where $\tilde\sigma_1:
\R^{n_1}\ra\R^m\otimes\R^{n_1}$ and $\tilde\sigma_2:
\R^{n}\ra\R^m\otimes\R^{n_2}$. And we also have a vector field $\tilde b=(\tilde b_1,
\tilde b_2)$ with the same structure of $b$ given above.

\begin{lemma}\label{sect-3-lem-2}
Suppose that $\sigma_1,\tilde\sigma_1\in W^{1,2q}_{x_1,loc}$ and $b_1,\tilde b_1\in W^{1,q}_{x_1,loc}$.
Moreover, $\sigma_2,\tilde\sigma_2\in W^{1,2q}_{x_1,x_2,loc}$ and $b_2,\tilde b_2\in W^{1,q}_{x_1,x_2,loc}$.
Let $X_t$ (resp. $\tilde X_t$) be the stochastic flow associated to the It\^o SDE \eqref{sect-3.1}
with coefficients $\sigma$ and $b$ (resp. $\tilde\sigma$ and $\tilde b$). Denote by
$\rho_t$ (resp. $\tilde \rho_t$) the Radon--Nikodym density of $X_t$ (resp. $\tilde X_t$)
with respect to $\mu$. Assume that
  $$\Lambda_{p,T}:=\sup_{0\leq t\leq T}\big(\|\rho_t\|_{L^p(\P\times\mu)}\vee
  \|\tilde \rho_t\|_{L^p(\P\times\mu)}\big)<+\infty,$$
where $p$ is the conjugate number of $q$. Then for any $\delta>0$,
  \begin{align*}
  &\E\int_{G_R\cap \tilde G_R}\log\bigg(\frac{\|X_{2,\cdot}-\tilde X_{2,\cdot}\|^2_{\infty,T}}{\delta^2}+1\bigg)\d\mu
  \leq \E\int_{G_R\cap \tilde G_R}\log\bigg(\frac{\|X-\tilde X\|^2_{\infty,T}}{\delta^2}+1\bigg)\d\mu\cr
  &\hskip6pt \leq C_T\Lambda_{p,T}\bigg\{C_{n,q}\Big[\|\nabla b\|_{L^q(B(3R))}
  +\|\nabla \sigma\|_{L^{2q}(B(3R))}+\|\nabla\sigma\|^2_{L^{2q}(B(3R))}\Big]\cr
  &\hskip60pt +\frac1{\delta^2}\|\sigma-\tilde\sigma\|^2_{L^{2q}(B(R))}
  +\frac1\delta \Big[\|\sigma-\tilde\sigma\|_{L^{2q}(B(R))}+
  \|b-\tilde b\|_{L^{q}(B(R))}\Big]\bigg\},
  \end{align*}
where $\tilde G_R:=\big\{(\omega,x)\in \Omega\times\R^n:\|\tilde X_\cdot(\omega,x)\|_{\infty,T}\leq R\big\}$
is the level set of the flow $\tilde X_t$.
\end{lemma}

We omit the proof here. The reader can consult \cite[Theorem 5.2]{FangLuoThalmaier}
for details. We mention that by Lemma \ref{sect-2-lem-2}, similar result for the
term $\|X_{1,\cdot}-\tilde X_{1,\cdot}\|_{\infty,T}$ holds. Next we focus on the existence
part of Theorem \ref{auxiliary-theorem} which needs to regularize the coefficients
$\sigma_1,b_1$ and $\sigma_2,b_2$ separately.

Let $\chi_1\in
C_c^\infty(\R^{n_1},\R_+)$ be such that $\int_{\R^{n_1}}\chi_1(x_1)\,\d x_1=1$ and
its support $\supp(\chi_1)\subset B_1(1)$, where $B_1(r)$ is a ball
in $\R^{n_1}$ centered at the origin with radius $r>0$. For $k\geq1$, define
$\chi_{1,k}(x_1)=k^{n_1}\chi_1(k x_1)$.
Next choose $\psi_1\in C_c^\infty(\R^{n_1},[0,1])$ so that $\psi_1|_{B_1(1)}
\equiv 1$ and $\psi_1$ vanishes outside $B_1(2)$. Denote by $\psi_{1,k}(x_1)=
\psi_1(x_1/k)$ for $k\geq1$. Now we set
  $$\sigma_{1,k}=(\sigma_1\ast\chi_{1,k})\,\psi_{1,k},\quad
  b_{1,k}=(b_1\ast\chi_{1,k})\,\psi_{1,k};$$
and
  \begin{equation}\label{sect-3.2}
  \sigma_{2,k}=(\sigma_2\ast\chi_k)\,\psi_k,\quad
  b_{2,k}=(b_2\ast\chi_k)\,\psi_k.
  \end{equation}
Here $\chi_k$ and $\psi_k$ are the same as in Section 2.
Then the coefficients $\sigma_{i,k}, b_{i,k}\in C_b^\infty(\R^n)$ ($i=1,2$).
Furthermore, by Lemma \ref{appendix-lem}, it holds
  \begin{equation}\label{auxiliary-theorem.2}
  \frac{|\sigma_{1,k}|}{1+|x_1|}\leq 2|\bar\sigma_1|\ast\chi_{1,k},
  \quad \frac{|b_{1,k}|}{1+|x_1|}\leq 2|\bar b_1|\ast\chi_{1,k}
  \end{equation}
and
  \begin{equation}\label{auxiliary-theorem.3}
  \frac{|\sigma_{2,k}|}{1+|x|}\leq 2|\bar\sigma_2|\ast\chi_k,
  \quad  \frac{|b_{2,k}|}{1+|x|}\leq 2|\bar b_2|\ast\chi_k.
  \end{equation}

We now consider the It\^o SDEs
  \begin{equation*}
  \begin{cases}
  \d X^k_{1,t}=\sigma_{1,k}(X^k_{1,t})\,\d B_t+ b_{1,k}(X^k_{1,t})\,\d t, & X^k_{1,0}=x_1,\\
  \d X^k_{2,t}=\sigma_{2,k}(X^k_{1,t},X^k_{2,t})\,\d B_t+ b_{2,k}(X^k_{1,t},X^k_{2,t})\,\d t, & X^k_{2,0}=x_2.
  \end{cases}
  \end{equation*}
For any $k\geq1$, the above equation determines a unique stochastic flow
$X^k_t=(X^k_{1,t},X^k_{2,t})$ of diffeomorphisms on $\R^n$. Moreover,
denoting by $\rho^k_t=\frac{\d[X^k_t]_\#\mu}{\d\mu}$,
then by Lemma \ref{density-estimate}, we have for any $p>1$ and $t\in [0,T]$,
  \begin{equation}\label{auxiliary-theorem.4}
  \|\rho^k_t\|_{L^p(\P\times\mu)}\leq \mu(\R^n)^{\frac1{p+1}}
  \bigg(\sup_{0\leq t\leq T}\int_{\R^n}\exp\big(tp^3|\Lambda_1^{\sigma_k}|^2
  -tp^2\Lambda_2^{\sigma_k,b_k}\big)\d\mu\bigg)^{\frac1{p(p+1)}},
  \end{equation}
where $\sigma_k=(\sigma_{1,k},\sigma_{2,k})$ and $b_k=(b_{1,k},b_{2,k})$.
We shall find a uniform estimate for the densities $\rho^k_t$, hence we need
the following lemma which is an analogue of Lemma \ref{sect-2-lem-3}.

\begin{lemma}\label{sect-3-lem-3}
There is a positive constant $C_0>0$ independent on $k\geq 1$, such that
  $$|\Lambda_1^{\sigma_k}|^2\leq C_0\big(|\div_{x_1}(\sigma_1)|^2
  +|\bar\sigma_1|^2\big)\ast\chi_{1,k}
  +C_0\big(|\div_{x_2}(\sigma_2)|^2+|\bar\sigma_2|^2\big)\ast\chi_k,$$
and
  \begin{align*}
  -\Lambda_2^{\sigma_k,b_k}&\leq C_0\big([\div_{x_1}(b_1)]^- +|\bar b_1|
  +|\bar\sigma_1|^2+|\nabla_{x_1}\sigma_1|^2\big)\ast\chi_{1,k}\cr
  &\quad +C_0\big([\div_{x_2}(b_2)]^- +|\bar b_2|
  +|\bar\sigma_2|^2+|\nabla_{x_2}\sigma_2|^2\big)\ast\chi_k.
  \end{align*}
\end{lemma}

\noindent{\bf Proof.} The proof is similar to Lemma \ref{sect-2-lem-3}. Indeed,
note that $\div(b_k)=\div_{x_1}(b_{1,k})+\div_{x_2}(b_{2,k})$ and we deal with
the two terms separately as in the proof of Lemma \ref{sect-2-lem-3}. The other
estimates can be established in the same way. Thanks to \eqref{key-observation},
the partial derivatives $\nabla_{x_1}\sigma_2$ do not show up here. \fin

\begin{lemma}[Uniform density estimate]\label{sect-3-lem-4}
For fixed $p>1$, there are two positive constants $C_{1,p}, C_{2,p}>0$ and $T_0>0$
small enough such that
  \begin{align*}
  \sup_{0\leq t\leq T_0}\|\rho^k_t\|_{L^p(\P\times\mu)}
  &\leq C_{1,p}\bigg(\int_{\R^{n_1}}\exp\big[C_{2,p}T_0\big([\div_{x_1}(b_1)]^- +|\bar b_1|
  +|\nabla_{x_1}\sigma_1|^2+|\bar\sigma_1|^2\big)\big]\d \mu_1\bigg)^{\frac1{p(p+1)}}\cr
  &\quad \times \bigg(\int_{\R^{n_1}}\exp\big[C_{2,p}T_0\big([\div_{x_2}(b_2)]^- +|\bar b_2|
  +|\nabla_{x_2}\sigma_2|^2+|\bar\sigma_2|^2\big)\big]\d\mu\bigg)^{\frac1{p(p+1)}}.
  \end{align*}
\end{lemma}

\noindent{\bf Proof.} Note that
$|\div_{x_i}(\sigma_i)|\leq |\nabla_{x_i}\sigma_i|\ (i=1,2)$, thus
the first estimate in Lemma \ref{sect-3-lem-3} becomes
  $$|\Lambda_1^{\sigma_k}|^2\leq C_0\big(|\nabla_{x_1}\sigma_1|^2
  +|\bar\sigma_1|^2\big)\ast\chi_{1,k}
  +C_0\big(|\nabla_{x_2}\sigma_2|^2+|\bar\sigma_2|^2\big)\ast\chi_k.$$
For any $t\in [0,T]$, the above inequality plus the second one in Lemma \ref{sect-3-lem-2} gives us
  \begin{align*}
  tp^3|\Lambda_1^{\sigma_k}|^2-tp^2\Lambda_2^{\sigma_k,b_k}
  &\leq 2Tp^3C_0\big([\div_{x_1}(b_1)]^- +|\bar b_1|+|\bar\sigma_1|^2
  +|\nabla_{x_1}\sigma_1|^2\big)\ast\chi_{1,k}\cr
  &\quad + 2Tp^3C_0\big([\div_{x_2}(b_2)]^- +|\bar b_2|+|\bar\sigma_2|^2
  +|\nabla_{x_2}\sigma_2|^2\big)\ast\chi_k.
  \end{align*}
Denote by
  $$\Phi_i=2Tp^3C_0\big([\div_{x_i}(b_i)]^- +|\bar b_i|+|\bar\sigma_i|^2
  +|\nabla_{x_i}\sigma_i|^2\big),\quad i=1,2.$$
Then $\Phi_1$ is a function defined on $\R^{n_1}$, while $\Phi_2$ is a function on the
whole $\R^n$. Now we have by Cauchy's inequality,
  \begin{align}\label{density-estimate.1.5}
  \hskip-20pt\int_{\R^n}\exp\big(tp^3|\Lambda_1^{\sigma_k}|^2
  -tp^2\Lambda_2^{\sigma_k,b_k}\big)\d\mu
  &\leq \int_{\R^n} e^{\Phi_1\ast\chi_{1,k}}\, e^{\Phi_2\ast\chi_k}\,\d\mu\cr
  \hskip-20pt&\leq \bigg[\int_{\R^n} e^{2\Phi_1\ast\chi_{1,k}}\,\d\mu \bigg]^{\frac12}
  \bigg[\int_{\R^n} e^{2\Phi_2\ast\chi_k}\,\d\mu \bigg]^{\frac12}.
  \end{align}

In the following we estimate the two integrals given in \eqref{density-estimate.1.5}.
First we have
  $$(1+|x|^2)^\alpha \geq (1+|x_1|^2)^{\alpha_1}\times(1+|x_2|^2)^{\alpha-\alpha_1}.$$
Thus
  \begin{align}\label{density-estimate.2}
  \int_{\R^n} e^{2\Phi_1\ast\chi_{1,k}}\,\d\mu
  &\leq \int_{\R^n}e^{2(\Phi_1\ast\chi_{1,k})(x_1)}\frac{\d x_1}{(1+|x_1|^2)^{\alpha_1}}
  \cdot\frac{\d x_2}{(1+|x_2|^2)^{\alpha-\alpha_1}}\cr
  &=\mu_2(\R^{n_2})\int_{\R^{n_1}}e^{2(\Phi_1\ast\chi_{1,k})(x_1)+\lambda_1(x_1)}\,\d x_1,
  \end{align}
where $\d\mu_2(x_2)=\frac{\d x_2}{(1+|x_2|^2)^{\alpha-\alpha_1}}$ is a finite
measure on $\R^{n_2}$ and $\lambda_1(x_1)=-\alpha_1\log(1+|x_1|^2)$.
Similar to \eqref{sect-2-lem-4.3}, we can show that there is a constant $C>0$
such that for any $k\geq 1$,
  \begin{equation}\label{density-estimate.3}
  \lambda_1(x_1)\leq (\lambda_1\ast\chi_{1,k})(x_1)+C\quad\mbox{for all } x_1\in\R^{n_1}.
  \end{equation}
Substituting \eqref{density-estimate.3} into the inequality \eqref{density-estimate.2}
and by Jensen's inequality, we obtain
  \begin{align}\label{density-estimate.4}
  \int_{\R^n} e^{2\Phi_1\ast\chi_{1,k}}\,\d\mu
  &\leq \mu_2(\R^{n_2})\, e^C\int_{\R^{n_1}}e^{[(2\Phi_1+\lambda_1)\ast\chi_{1,k}](x_1)}\,\d x_1\cr
  &\leq \mu_2(\R^{n_2})\, e^C\int_{\R^{n_1}}\big[(e^{2\Phi_1+\lambda_1})\ast\chi_{1,k}\big](x_1)\,\d x_1\cr
  &=\mu_2(\R^{n_2})\, e^C \int_{\R^{n_1}}e^{2\Phi_1}\,\d\mu_1.
  \end{align}

The second integral on the right hand side of \eqref{density-estimate.1.5} can be treated in
a similar way, thanks to \eqref{sect-2-lem-4.3}. Hence
  \begin{equation}\label{density-estimate.5}
  \int_{\R^n} e^{2\Phi_2\ast\chi_k}\,\d\mu
  \leq e^{\bar C}\int_{\R^{n}}e^{2\Phi_2}\,\d\mu.
  \end{equation}
Now combining the inequalities \eqref{density-estimate.1.5}, \eqref{density-estimate.4}
and \eqref{density-estimate.5}, we finally obtain from the
definition of $\Phi_1$ and $\Phi_2$ that
  \begin{align*}
  &\int_{\R^n}\exp\big(tp^3|\Lambda_1^{\sigma_k}|^2
  -tp^2\Lambda_2^{\sigma_k,b_k}\big)\d\mu\cr
  &\quad \leq \big(\mu_2(\R^{n_2})\, e^{C+\bar C}\big)^{\frac12}
  \bigg[\int_{\R^{n_1}}\exp\Big\{4Tp^3C_0\big([\div_{x_1}(b_1)]^- +|\bar b_1|+|\bar\sigma_1|^2
  +|\nabla_{x_1}\sigma_1|^2\big)\Big\}\,\d\mu_1\bigg]^{\frac1{2}}\cr
  &\quad\quad \times \bigg[\int_{\R^{n}}\exp\Big\{4Tp^3C_0\big([\div_{x_2}(b_2)]^- +|\bar b_2|+|\bar\sigma_2|^2
  +|\nabla_{x_2}\sigma_2|^2\big)\Big\}\,\d\mu\bigg]^{\frac1{2}}.
  \end{align*}
Substituting this inequality into \eqref{auxiliary-theorem.4}, we see that for any
$k\geq 1$,
  \begin{align*}
  \sup_{t\leq T}\|\rho^k_t\|_{L^p(\P\times\mu)}&\leq
  C_{1,p}\bigg[\int_{\R^{n_1}}\exp\Big\{C_{2,p}T\big([\div_{x_1}(b_1)]^- +|\bar b_1|+|\bar\sigma_1|^2
  +|\nabla_{x_1}\sigma_1|^2\big)\Big\}\,\d\mu_1\bigg]^{\frac1{2p(p+1)}}\cr
  &\quad \times \bigg[\int_{\R^n}\exp\Big\{C_{2,p}T\big([\div_{x_2}(b_2)]^- +|\bar b_2|+|\bar\sigma_2|^2
  +|\nabla_{x_2}\sigma_2|^2\big)\Big\}\,\d\mu\bigg]^{\frac1{2p(p+1)}},
  \end{align*}
where $C_{1,p},C_{2,p}$ are two positive constants independent on $k$ and $T$.
Under the conditions (H2) and (H4), there exists $T_0>0$ small enough such
that the quantity on the right hand side is finite. \fin

\medskip

Having Lemma \ref{sect-3-lem-4} in hand, we can follow the line of
arguments in Section 2 to prove Theorem \ref{auxiliary-theorem}. We omit the
details.

\section{SDE with partially Sobolev coefficients}

In this section we aim at generalizing Theorem \ref{auxiliary-theorem}
to the case where the coefficients $\sigma_2$ and $b_2$ only have partial
Sobolev regularity. More precisely, we replace the condition (H3) by
  \begin{itemize}
  \item[\rm(H3$^\prime$)] $\sigma_2\in L^{2q}_{x_1,loc}(W^{1,2q}_{x_2,loc}),\,
  b_2\in L^q_{x_1,loc}(W^{1,q}_{x_2,loc})$,
  \end{itemize}
and we shall show that the results of Theorem \ref{auxiliary-theorem}
still hold.

To achieve such an extension, we need an a-priori estimate which is analogous
to Lemma \ref{sect-3-lem-2}, but only involving partial derivatives of $\sigma_2$
and $b_2$. First we introduce some notations.
Throughout this section we fix a pair of functions
  $$\sigma_1:\R^{n_1}\ra \R^m\otimes\R^{n_1}\quad \mbox{and}\quad
  b_1:\R^{n_1}\ra \R^{n_1}$$
which satisfy the assumptions (H1) and (H2) in Section 3. Under these
conditions, it is known that the following It\^o SDE
  $$\d X_{1,t}=\sigma_{1}(X_{1,t})\,\d B_t+ b_1(X_{1,t})\,\d t,\quad X_{1,0}=x_1$$
generates a unique stochastic flow of measurable maps on $\R^{n_1}$, which
leaves the reference measure $\mu_1$ absolutely continuous, as shown in Theorem \ref{Extension}.

Let
  $$\sigma_2,\tilde\sigma_2:\R^n \ra \R^m\otimes\R^{n_2}\quad \mbox{and} \quad
  b_2,\tilde b_2:\R^n\ra \R^{n_2}$$
be measurable functions, all verifying the conditions (H3$^\prime$). Denote by
  $$\sigma=(\sigma_1,\sigma_2),\
  b=(b_1,b_2)
  \quad \mbox{and} \quad
  \tilde\sigma=(\sigma_1,\tilde\sigma_2),\
  \tilde b=(b_1,\tilde b_2).$$
Let $X_t=(X_{1,t},X_{2,t})$ (resp. $\tilde X_t=(X_{1,t},\tilde X_{2,t})$) be the
stochastic flow generated by the It\^o SDE
\eqref{Ito-SDE} with coefficients $\sigma$ and $b$ (resp. $\tilde\sigma$ and $\tilde b$).

\begin{lemma}[A-priori estimate]\label{a-priori-estimate}
Suppose that for any $t\in[0,T]$, the push-forwards $(X_t)_\#\mu$ and $(\tilde X_t)_\#\mu$ of the
reference measure $\mu$ are absolutely continuous with respect to itself, with
density functions $\rho_t$ and $\tilde\rho_t$ respectively. Moreover,
  \begin{equation}\label{a-priori-estimate.0}
  \Lambda_{p,T}:=\sup_{0\leq t\leq T}\|\rho_t\|_{L^p(\P\otimes\mu)}\vee
  \|\tilde\rho_t\|_{L^p(\P\otimes\mu)}<+\infty,
  \end{equation}
where $p$ is the conjugate number of $q$. Then for any $\delta>0$,
  \begin{align*}
  &\E\int_{G_R\cap \tilde G_R}\log\bigg(\frac{\|X_{2}-\tilde X_{2}\|^2_{\infty,T}}{\delta^2}+1\bigg)\d\mu\cr
  &\hskip6pt \leq C_T\Lambda_{p,T}\bigg\{C_{n_2,q}\Big[\|\nabla_{x_2}b_2\|_{L^q(B(4R))}
  +\|\nabla_{x_2}\sigma_2\|_{L^{2q}(B(4R))}+\|\nabla_{x_2}\sigma_2\|^2_{L^{2q}(B(4R))}\Big]\cr
  &\hskip60pt +\frac1{\delta^2}\|\sigma_2-\tilde\sigma_2\|^2_{L^{2q}(B(R))}
  +\frac1\delta \Big[\|\sigma_2-\tilde\sigma_2\|_{L^{2q}(B(R))}+
  \|b_2-\tilde b_2\|_{L^{q}(B(R))}\Big]\bigg\},
  \end{align*}
where $G_R$ and $\tilde G_R$ are the level sets of $X_t$ and $\tilde X_t$ respectively.
\end{lemma}

\noindent{\bf Proof.} We follow the idea of the proof of \cite[Theorem 5.2]{FangLuoThalmaier}
(see also \cite[Lemma 4.1]{Zhang12}). Denote by $\xi_t=X_{2,t}-\tilde X_{2,t}$.
Then $\xi_0=0$. By the It\^o formula,
  \begin{align}\label{a-priori-estimate.1}
  \d\log(|\xi_t|^2+\delta^2)&=\frac{2\big\<\xi_t,[\sigma_2(X_t)-\tilde \sigma_2(\tilde X_t)]\,\d B_t\big\>}
  {|\xi_t|^2+\delta^2} +\frac{2\<\xi_t,b_2(X_t)-\tilde b_2(\tilde X_t)\>}{|\xi_t|^2+\delta^2}\,\d t\cr
  &\hskip12pt +\frac{|\sigma_2(X_t)-\tilde \sigma_2(\tilde X_t)|^2}{|\xi_t|^2+\delta^2}\,\d t
  -\frac{2\big|[\sigma_2(X_t)-\tilde \sigma_2(\tilde X_t)]^\ast\xi_t\big|^2}{(|\xi_t|^2+\delta^2)^2}\,\d t\cr
  &=: \sum_{i=1}^4 \d I_i(t).
  \end{align}
Note that the last term is negative, hence we omit it. We shall estimate the
other terms in the sequel.

Let $\tau_R(x)=\inf\{t\geq0:|X_t(x)|\vee|\tilde X_t(x)|>R\}$ for $x\in\R^n$.
Remark that almost surely, $G_R,\tilde G_R\subset\{x:\tau_R(x)>T\}$ and for any $t\geq0$,
$\{\tau_R>t\}\subset B(R)$. Thus by Cauchy's inequality,
  \begin{align*}
  \E\bigg[\int_{G_R\cap\tilde G_R}\sup_{0\leq t\leq T}|I_1(t)|\,\d\mu\bigg]
  &\leq \E\bigg[\int_{B(R)}\sup_{0\leq t\leq T\wedge\tau_R}|I_1(t)|\,\d\mu\bigg]\cr
  &\leq \mu(\R^n)^{\frac12}\bigg[\int_{B(R)}\E\bigg(\sup_{0\leq t\leq T\wedge\tau_R}|I_1(t)|^2\bigg)\d\mu\bigg]^{\frac12}.
  \end{align*}
Burkholder's inequality gives us
  \begin{align*}
  \E\bigg(\sup_{0\leq t\leq T\wedge\tau_R}|I_1(t)|^2\bigg)
  &\leq 16\,\E\bigg(\int_0^{T\wedge\tau_R}\frac{\big|[\sigma_2(X_t)-\tilde \sigma_2(\tilde X_t)]^\ast\xi_t\big|^2}
  {(|\xi_t|^2+\delta^2)^2}\,\d t\bigg)\cr
  &\leq 16\,\E\bigg(\int_0^{T\wedge\tau_R}\frac{|\sigma_2(X_t)-\tilde \sigma_2(\tilde X_t)|^2}
  {|\xi_t|^2+\delta^2}\,\d t\bigg).
  \end{align*}
As a result, by changing the order of integration, we obtain
  \begin{align}\label{a-priori-estimate.2}
  \E\bigg[\int_{G_R\cap\tilde G_R}\sup_{0\leq t\leq T}|I_1(t)|\,\d\mu\bigg]
  &\leq 4C_{\alpha,n} \bigg[\int_{B(R)}\E\bigg(\int_0^{T\wedge\tau_R}\frac{|\sigma_2(X_t)-\tilde \sigma_2(\tilde X_t)|^2}
  {|\xi_t|^2+\delta^2}\,\d t\bigg)\d\mu\bigg]^{\frac12}\cr
  &=4C_{\alpha,n} \bigg[\int_0^{T}\bigg(\E\int_{\{\tau_R>t\}}\frac{|\sigma_2(X_t)-\tilde \sigma_2(\tilde X_t)|^2}
  {|\xi_t|^2+\delta^2}\,\d\mu\bigg)\d t\bigg]^{\frac12}.
  \end{align}
Note that
  $$\sigma_2(X_t)-\tilde \sigma_2(\tilde X_t)
  =\sigma_2(X_{t})-\sigma_2(\tilde X_{t})
  +\sigma_2(\tilde X_{t})-\tilde\sigma_2(\tilde X_{t}).$$
We have by \eqref{a-priori-estimate.0} and H\"older's inequality that
  \begin{align*}
  \E\int_{\{\tau_R>t\}}\frac{|\sigma_2(\tilde X_t)-\tilde \sigma_2(\tilde X_t)|^2}
  {|\xi_t|^2+\delta^2}\,\d\mu
  &\leq \frac1{\delta^2}\,\E\int_{\{\tau_R>t\}}
  \big|(\sigma_2-\tilde\sigma_2)\ch_{B(R)}\big|^2(\tilde X_t)\,\d\mu\cr
  &\leq \frac1{\delta^2}\,\E\int_{B(R)}|\sigma_2-\tilde\sigma_2|^2\tilde \rho_t\,\d\mu\cr
  &\leq \frac{\Lambda_{p,T}}{\delta^2}\|\sigma_2-\tilde\sigma_2\|^2_{L^{2q}(B(R),\mu)}.
  \end{align*}
Since $\mu|_{B(R)}\leq \L_n|_{B(R)}$ for any $R>0$, we obtain
  \begin{equation}\label{a-priori-estimate.3}
  \E\int_{\{\tau_R>t\}}\frac{|\sigma_2(\tilde X_t)-\tilde \sigma_2(\tilde X_t)|^2}
  {|\xi_t|^2+\delta^2}\,\d\mu
  \leq \frac{\Lambda_{p,T}}{\delta^2}\|\sigma_2-\tilde\sigma_2\|^2_{L^{2q}(B(R))}.
  \end{equation}

Next on the set $\{\tau_R>t\}$, we have $X_t,\tilde X_t\in B(R)$, hence $|X_t-\tilde X_t|_{\R^n}
=|X_{2,t}-\tilde X_{2,t}|_{\R^{n_2}}\leq 2R$. As $(X_t)_\#\mu\ll \mu$ and
$(\tilde X_t)_\#\mu\ll \mu$, we can apply Lemma \ref{pointwise-inequality}(i) to get
  $$|\sigma_2(X_{t})-\sigma_2(\tilde X_{t})|\leq C_{n_2}|X_{2,t}-\tilde X_{2,t}|
  \,\big(M_{2,2R}|\nabla_{x_2}\sigma_2|(X_t)+M_{2,2R}|\nabla_{x_2}\sigma_2|(\tilde X_t)\big).$$
Thus
  \begin{align*}
  \E\int_{\{\tau_R>t\}}\frac{|\sigma_2(X_t)-\sigma_2(\tilde X_t)|^2}
  {|\xi_t|^2+\delta^2}\,\d\mu
  &\leq C^2_{n_2}\E\int_{\{\tau_R>t\}}\big(M_{2,2R}|\nabla_{x_2}\sigma_2|(X_t)
  +M_{2,2R}|\nabla_{x_2}\sigma_2|(\tilde X_t)\big)^2\,\d\mu\cr
  &\leq 2C^2_{n_2}\E\int_{B(R)}\big(M_{2,2R}|\nabla_{x_2}\sigma_2|\big)^2(\rho_t+\tilde\rho_t)\,\d\mu.
  \end{align*}
H\"older's inequality gives us
  \begin{align}\label{a-priori-estimate.4}
  \E\int_{\{\tau_R>t\}}\frac{|\sigma_2(X_t)-\sigma_2(\tilde X_t)|^2}
  {|\xi_t|^2+\delta^2}\,\d\mu
  &\leq 4C^2_{n_2}\Lambda_{p,T}\bigg(\int_{B(R)}
  \big(M_{2,2R}|\nabla_{x_2}\sigma_2|\big)^{2q}\,\d\mu\bigg)^{\frac1q}.
  \end{align}
We have
  \begin{align*}
  \int_{B(R)}\big(M_{2,2R}|\nabla_{x_2}\sigma_2|\big)^{2q}\,\d\mu
  &\leq \int_{B(R)}\big(M_{2,2R}|\nabla_{x_2}\sigma_2|\big)^{2q}\,\d x\cr
  &\leq \int_{B_1(R)}\d x_1\!\int_{B_2(R)}\big(M_{2,2R}|\nabla_{x_2}\sigma_2|\big)^{2q}\,\d x_2.
  \end{align*}
Recall that $B_i(R)$ is a ball in $\R^{n_i}$ centered at the origin with radius $R$, $i=1,2$.
Lemma \ref{pointwise-inequality}(ii) gives us
  $$\int_{B_2(R)}\big(M_{2,2R}|\nabla_{x_2}\sigma_2|\big)^{2q}\,\d x_2
  \leq C_{q,n_2}\int_{B_2(3R)}|\nabla_{x_2}\sigma_2|^{2q}\,\d x_2.$$
Therefore
  \begin{align*}
  \int_{B(R)}\big(M_{2,2R}|\nabla_{x_2}\sigma_2|\big)^{2q}\,\d\mu
  \leq C_{q,n_2} \int_{B(4R)}|\nabla_{x_2}\sigma_2|^{2q}\,\d x.
  \end{align*}
Substituting this estimate into \eqref{a-priori-estimate.4}, we obtain
  \begin{align*}
  \E\int_{\{\tau_R>t\}}\frac{|\sigma_2(X_t)-\sigma_2(\tilde X_t)|^2}
  {|\xi_t|^2+\delta^2}\,\d \mu
  & \leq C'_{q,n_2}\Lambda_{p,T}\bigg(\int_{B(4R)}
  |\nabla_{x_2}\sigma_2|^{2q}\,\d x\bigg)^{\frac1q}\cr
  &=C'_{q,n_2}\Lambda_{p,T}\|\nabla_{x_2}\sigma_2\|^2_{L^{2q}(B(4R))}.
  \end{align*}
Combining this inequality with \eqref{a-priori-estimate.2} and
\eqref{a-priori-estimate.3}, we arrive at
  \begin{align}\label{a-priori-estimate.6}
  &\E\bigg[\int_{G_R\cap\tilde G_R}\sup_{0\leq t\leq T}|I_1(t)|\,\d\mu\bigg]\cr
  &\quad \leq C_T\Lambda_{p,T}^{\frac12}\bigg[\frac1{\delta^2}\|\sigma_2-\tilde\sigma_2\|^2_{L^{2q}(B(R))}
  +C'_{q,n_2}\|\nabla_{x_2}\sigma_2\|^2_{L^{2q}(B(4R))}\bigg]^{\frac12}.
  \end{align}

Now we begin estimating the term $I_2(t)$. We have
  $$\E\bigg[\int_{G_R\cap\tilde G_R}\sup_{0\leq t\leq T}|I_2(t)|\,\d\mu\bigg]
  \leq 2\int_0^T\bigg[\E\int_{G_R\cap\tilde G_R}\frac{|b_2(X_t)-\tilde b_2(\tilde X_t)|}
  {(|\xi_{2,t}|^2+\delta^2)^{\frac12}}\,\d\mu\bigg]\d t.$$
For $x\in G_R\cap\tilde G_R$, one has $\tilde X_t(x)\in B(R)$ for all $t\in[0,T]$, then
  \begin{equation}\label{a-priori-estimate.7}
  \E\int_{G_R\cap\tilde G_R}\frac{|b_2(\tilde X_t)-\tilde b_2(\tilde X_t)|}
  {(|\xi_{2,t}|^2+\delta^2)^{\frac12}}\,\d\mu
  \leq \frac1\delta \,\E\int_{B(R)}|b_2-\tilde b_2|\tilde\rho_t\,\d\mu
  \leq \frac{\Lambda_{p,T}}{\delta}\|b_2-\tilde b_2\|_{L^q(B(R))}.
  \end{equation}
By Lemma \ref{pointwise-inequality}(i) and H\"older's inequality,
analogous arguments as for estimating \eqref{a-priori-estimate.4} leads to
  \begin{align*}
  \E\int_{G_R\cap\tilde G_R}\frac{|b_2(X_t)-b_2(\tilde X_t)|}
  {(|\xi_{2,t}|^2+\delta^2)^{\frac12}}\,\d\mu
  &\leq C_{n_2}\E\int_{G_R\cap\tilde G_R}\big(M_{2,2R}|\nabla_{x_2}b_2|(X_t)
  +M_{2,2R}|\nabla_{x_2}b_2|(\tilde X_t)\big)\d\mu\cr
  &\leq C_{n_2}\E\int_{B(R)}\big(M_{2,2R}|\nabla_{x_2}b_2|\big)(\rho_t+\tilde\rho_t)\,\d\mu\cr
  &\leq 2C''_{q,n_2}\Lambda_{p,T}\|\nabla_{x_2}b_2\|_{L^q(B(4R))}.
  \end{align*}
This together with \eqref{a-priori-estimate.7} gives us
  \begin{equation}\label{a-priori-estimate.8}
  \E\bigg[\int_{G_R\cap\tilde G_R}\sup_{0\leq t\leq T}|I_2(t)|\,\d\mu\bigg]
  \leq 2T\Lambda_{p,T}\bigg(\frac1\delta \|b_2-\tilde b_2\|_{L^q(B(R))}
  +C''_{q,n_2}\|\nabla_{x_2}b_2\|_{L^q(B(4R))}\bigg).
  \end{equation}
Similarly we can show that
  \begin{equation}\label{a-priori-estimate.9}
  \E\bigg[\int_{G_R\cap\tilde G_R}\sup_{0\leq t\leq T}|I_3(t)|\,\d\mu\bigg]
  \leq CT\Lambda_{p,T}\bigg(\frac1{\delta^2}\|\sigma_2-\tilde\sigma_2\|^2_{L^{2q}(B(R))}
  +C'_{q,n_2}\|\nabla_{x_2}\sigma_2\|^2_{L^{2q}(B(4R))}\bigg).
  \end{equation}

Combining the estimates \eqref{a-priori-estimate.6}, \eqref{a-priori-estimate.8}
and  \eqref{a-priori-estimate.9}, we obtain the result. \fin

\medskip

The a-priori estimate in Lemma \ref{a-priori-estimate} has some direct consequences.
The first one is the stability of generalized stochastic flow,
which is the content of the following theorem.

\begin{theorem}[Stability]\label{sect-4-stability}
Suppose there is a sequence of coefficients $\sigma_{2,k}:\R^n\ra \R^m\otimes\R^{n_2}$
and $b_{2,k}:\R^n\ra \R^{n_2}$, verifying the conditions (H3$'$) and (H4).
Assume that $\sigma_{2,k}$ (resp. $b_{2,k}$) converge to $\sigma_2$ (resp. $b_2$) in
$L^{2q}_{loc}(\R^n)$ (resp. $L^q_{loc}(\R^n)$) as $k\ra\infty$. We also assume that
  \begin{align}\label{sect-4-stability.1}
  C_1:=\sup_{k\geq1}\big[\|\sigma_{2,k}\|_{L^{2q}(\mu)}+\|b_{2,k}\|_{L^{q}(\mu)}\big]<+\infty,
  \end{align}
and for any $R>0$,
  \begin{align}\label{sect-4-stability.2}
  C_{2,R}:=\sup_{k\geq1}\big[\|\nabla_{x_2}b_{2,k}\|_{L^q(B(R))}
  +\|\nabla_{x_2}\sigma_{2,k}\|_{L^{2q}(B(R))}\big]<+\infty.
  \end{align}

Let $X^k_t=(X_{1,t},X^k_{2,t})$ be the stochastic flow generated by the It\^o SDE \eqref{Ito-SDE}
with the coefficients $\sigma_k=(\sigma_1,\sigma_{2,k})$ and $b_k=(b_1,b_{2,k})$.
Suppose that for all $k\geq1$, the density function $\rho^k_t:=\frac{\d (X^k_t)_\#\mu}{\d\mu}$ exists and
  \begin{equation}\label{sect-4-stability.3}
  \Lambda_{p,T}:=\sup_{k\geq1}\sup_{0\leq t\leq T}\|\rho^k_t\|_{L^p(\P\times\mu)}<+\infty.
  \end{equation}
Then there exists a random field $X_2:\Omega\times\R^n\ra C([0,T],\R^{n_2})$ such that
  $$\lim_{k\ra\infty}\E\int_{\R^n}1\wedge \|X^k_2-X_2\|_{\infty,T}\,\d\mu=0.$$
\end{theorem}

\noindent{\bf Proof.} The proof is similar to that of Proposition \ref{sect-2-prop-1}.
For any $k\geq1$, let $G_R^k$ be the level set of the flow $X^k_t$:
  $$G_R^k=\big\{(\omega,x)\in \Omega\times\R^n:\|X^k_\cdot(\omega,x)\|_{\infty,T}\leq R\big\}. $$
Under the conditions \eqref{sect-4-stability.1} and \eqref{sect-4-stability.3}, we can apply
Lemma \ref{sect-3-lem-1} to get that
  $$(\P\times\mu)\big[(G_R^{k})^c\big] \leq \frac{C'_1}R,$$
where $C'_1$ depends only on $C_1$ and $\Lambda_{p,T}$.
As a result, for any $k,l\geq1$ and $R>0$,
  \begin{align}\label{sect-4-stability.4}
  (\P\times\mu)\big[(G_R^k\cap G_R^l)^c\big]&\leq
  (\P\times\mu)\big[(G_R^k)^c\big]+(\P\times\mu)\big[(G_R^l)^c\big]
  \leq \frac{2C'_1}R.
  \end{align}

Now applying Lemma \ref{a-priori-estimate} to the flows $X^k_t$ and $X^l_t$, we get
  \begin{align}\label{sect-4-stability.7}
  &\E\int_{G_R^k\cap G_R^l}\log\bigg(\frac{\|X^k_2- X^l_2\|^2_{\infty,T}}{\delta^2}+1\bigg)\d\mu\cr
  &\quad \leq C_T\Lambda_{p,T}\bigg\{C_{q,n_2}\Big[\|\nabla_{x_2}b_{2,k}\|_{L^q(B(4R))}
  +\|\nabla_{x_2}\sigma_{2,k}\|_{L^{2q}(B(4R))}+\|\nabla_{x_2}\sigma_{2,k}\|^2_{L^{2q}(B(4R))}\Big]\cr
  &\hskip22pt +\frac1{\delta^2}\|\sigma_{2,k}-\sigma_{2,l}\|^2_{L^{2q}(B(R))}
  +\frac1\delta \Big[\|\sigma_{2,k}-\sigma_{2,l}\|_{L^{2q}(B(R))}+
  \|b_{2,k}-b_{2,l}\|_{L^{q}(B(R))}\Big]\bigg\}.
  \end{align}
Since $\sigma_{2,k}\ra\sigma_2$ in $L^{2q}_{loc}(\R^n)$ and $b_{2,k}\ra b_2$ in $L^q_{loc}(\R^n)$
as $k\ra\infty$, we see that
  $$\delta_{k,l}:=\|\sigma_{2,k}-\sigma_{2,l}\|_{L^{2q}(B(R))}+\|b_{2,k}-b_{2,l}\|_{L^{q}(B(R))}\ra0$$
as $k,l$ goes to $\infty$. Next by \eqref{sect-4-stability.2}, there is a constant
$C'_R>0$ such that for all $k\geq 1$,
  $$\|\nabla_{x_2}b_{2,k}\|_{L^q(B(4R))}+\|\nabla_{x_2}\sigma_{2,k}\|_{L^{2q}(B(4R))}
  +\|\nabla_{x_2}\sigma_{2,k}\|^2_{L^{2q}(B(4R))}\leq C'_R.$$
Consequently, by taking $\delta=\delta_{k,l}$ in \eqref{sect-4-stability.4}, we
can find a positive constant $\hat C_0>0$ such that
  \begin{equation}\label{sect-4-stability.8}
  \E\int_{G_R^k\cap G_R^l}\log\bigg(\frac{\|X^k_2- X^l_2\|^2_{\infty,T}}{\delta_{k,l}^2}+1\bigg)\d\mu\leq\hat C_0
  \quad \mbox{for all } k,l\geq1.
  \end{equation}

Let $\eta\in (0,1)$ and define
  $$\Sigma_\eta^{k,l}=\{(\omega,x)\in \Omega\times\R^n:\|X^k_2- X^l_2\|_{\infty,T}\leq \eta\}.$$
Then
  \begin{align}\label{sect-4-stability.9}
  &\int_{G_R^k\cap G_R^l}1\wedge \|X^k_2-X^l_2\|_{\infty,T}\,\d(\P\times\mu)\cr
  &\quad =\bigg[\int_{(G_R^k\cap G_R^l)\cap\Sigma_\eta^{k,l}}+\int_{(G_R^k\cap G_R^l)\setminus\Sigma_\eta^{k,l}}\bigg]
  1\wedge \|X^k_2-X^l_2\|_{\infty,T}\,\d(\P\times\mu)\cr
  &\quad =:J_1+J_2.
  \end{align}
By Chebyshev's inequality and \eqref{sect-4-stability.8},
we have
  \begin{align*}
  J_2 \leq \frac1{\log\big(\frac{\eta^2}{\delta^2_{k,l}}+1\big)}
  \int_{G_R^k\cap G_R^l}\log\bigg(\frac{\|X^k_2- X^l_2\|^2_{\infty,T}}{\delta_{k,l}^2}+1\bigg)\d\mu
  \leq \frac{\hat C_0}{\log\big(\frac{\eta^2}{\delta^2_{k,l}}+1\big)}.
  \end{align*}
Therefore
  \begin{equation}\label{sect-4-stability.10}
  \limsup_{k,l\ra\infty}J_2=0.
  \end{equation}
On the other hand, $J_1\leq \eta\, \mu(\R^n)$. Combining this with \eqref{sect-4-stability.9}
and \eqref{sect-4-stability.10}, we obtain, by first letting $k,l\ra\infty$
and then $\eta\da0$, that
  $$\lim_{k,l\ra\infty}\int_{G_R^k\cap G_R^l}1\wedge \|X^k_2-X^l_2\|_{\infty,T}\,\d(\P\times\mu)=0$$
for any $R>0$. This together with \eqref{sect-4-stability.4} leads to the desired result. \fin

\medskip

Now we are ready to show the existence of generalized stochastic flows
to the It\^o SDE \eqref{Ito-SDE}.

\begin{theorem}[Existence]\label{sect-4-existence}
Under the assumptions (H1), (H2), (H3$'$) and (H4), the It\^o SDE \eqref{Ito-SDE} generates
a stochastic flow $X_t=(X_{1,t},X_{2,t})$, which is well defined on some small interval $[0,T_1]$.
Moreover, the Radon--Nikodym density $\rho_t:=\frac{\d (X_t)_\#\mu}{\d\mu}$ exists
and satisfies
  $$\sup_{0\leq t\leq T_1}\|\rho_t\|_{L^p(\P\times\mu)}<+\infty.$$
\end{theorem}

\noindent{\bf Proof.} We split the proof into three steps.

{\it Step 1}. In this step we shall regularize the coefficients
$\sigma_2,b_2$, and then apply Theorem \ref{auxiliary-theorem} to get a sequence of
stochastic flows.

To this end, we define $\sigma_{2,k}$ and $b_{2,k}$ as in \eqref{sect-3.2}.
We remark that there is no need to regularize the coefficients
$\sigma_1$ and $b_1$. Consider the family of It\^o's SDE:
  \begin{equation}\label{sect-4-existence.1}
  \begin{cases}
  \d X_{1,t}=\sigma_{1}(X_{1,t})\,\d B_t+ b_{1}(X_{1,t})\,\d t, & X_{1,0}=x_1,\\
  \d X^k_{2,t}=\sigma_{2,k}(X_{1,t},X^k_{2,t})\,\d B_t+ b_{2,k}(X_{1,t},X^k_{2,t})\,\d t, & X_{2,0}=x_2.
  \end{cases}
  \end{equation}

Now we check that the regularized coefficients  $\sigma_{2,k}$ and $b_{2,k}$
satisfy the conditions (H3) and (H4) stated at the beginning of Section 3.
Under the assumption (H3$'$), it is clear that
$\sigma_{2,k}\in W^{1,2q}_{x_1,x_2,loc}, b_{2,k}\in W^{1,q}_{x_1,x_2,loc}$,
hence (H3) is verified. Now we show that there is $p_1>0$ small enough such that
  \begin{align*}
  \int_{\R^{n}}\exp\big\{p_1\big([\div_{x_2}(b_{2,k})]^- +|\bar b_{2,k}|+|\bar\sigma_{2,k}|^2
  +|\nabla_{x_2}\sigma_{2,k}|^2\big)\big\}\,\d\mu<+\infty,
  \end{align*}
where $\bar b_{2,k}=\frac{b_{2,k}}{1+|x|}$ and $\bar \sigma_{2,k}=\frac{\sigma_{2,k}}{1+|x|}$.
In fact, similar to \eqref{sect-2-lem-3.4} and \eqref{sect-2-lem-3.6}, we have
  $$[\div_{x_2}(b_{2,k})]^- \leq \big([\div_{x_2}(b_2)]^- +2C|\bar b_2|\big)\ast\chi_k$$
and
  $$|\nabla_{x_2}\sigma_{2,k}|^2\leq C\big(|\nabla_{x_2}\sigma_2|^2+|\bar\sigma_2|^2\big)\ast\chi_k.$$
These estimates together with the inequalities \eqref{auxiliary-theorem.3} give us
  \begin{align}\label{sect-4-existence.2}
  &[\div_{x_2}(b_{2,k})]^- +|\bar b_{2,k}|+|\bar\sigma_{2,k}|^2+|\nabla_{x_2}\sigma_{2,k}|^2\cr
  &\quad \leq 2C\big([\div_{x_2}(b_2)]^- +|\bar b_2|+|\bar\sigma_2|^2+|\nabla_{x_2}\sigma_2|^2\big)\ast\chi_k.
  \end{align}
Now similar to the proof of Lemma \ref{sect-2-lem-4}, we can show that
  \begin{align}\label{sect-4-existence.2.5}
  &\int_{\R^n}\exp\big\{p\big([\div_{x_2}(b_{2,k})]^- +|\bar b_{2,k}|+|\bar\sigma_{2,k}|^2
  +|\nabla_{x_2}\sigma_{2,k}|^2\big)\big\}\,\d\mu\cr
  &\quad\leq \int_{\R^n}\exp\big\{2pC\big([\div_{x_2}(b_2)]^- +|\bar b_2|
  +|\bar\sigma_2|^2+|\nabla_{x_2}\sigma_2|^2\big)
  \ast\chi_k\big\}\,\d\mu\cr
  &\quad \leq 3^\alpha\int_{\R^n}\exp\big\{2pC\big([\div_{x_2}(b_2)]^- +|\bar b_2|
  +|\bar\sigma_2|^2+|\nabla_{x_2}\sigma_2|^2\big)\big\}\,\d\mu,
  \end{align}
where $C>0$ is independent on $k\geq 1$. Hence when $p\leq p_1:=p_0/2C$,
the right hand side is finite; in other words,
the condition (H4) is also satisfied.

Next, since $\sigma_1$ and $b_1$
satisfy (H1) and (H2), we can apply Theorem \ref{auxiliary-theorem} to conclude that
for every $k\geq1$, the It\^o SDE \eqref{sect-4-existence.1} generates a unique stochastic
flow $X^k_t=(X_{1,t},X^k_{2,t})$ which leaves the reference measure $\mu$ absolutely continuous,
and by Lemma \ref{sect-3-lem-3}, there is $T_0$ small enough such that the
Radon--Nikodym density $\rho^k_t:=\frac{\d(X^k_t)_\#\mu}{\d\mu}$ has the following estimate:
for all $t\leq T_0$,
  \begin{align*}
  \|\rho^k_t\|_{L^p(\P\times\mu)}&\leq
  C_{1,p}\bigg[\int_{\R^{n_1}} \exp\big\{C_{2,p}T_0\big([\div_{x_1}(b_1)]^-
  +|\bar b_1|+|\bar\sigma_1|^2+|\nabla_{x_1}\sigma_1|^2\big)\big\}\d\mu_1\bigg]^{\frac1{2p(p+1)}}\cr
  &\hskip12pt \times\bigg[\int_{\R^{n}} \exp\big\{ C_{2,p}T_0\big([\div_{x_2}(b_{2,k})]^-
  +|\bar b_{2,k}|+|\bar\sigma_{2,k}|^2+|\nabla_{x_2}\sigma_{2,k}|^2\big)\big\}\d\mu \bigg]^{\frac1{2p(p+1)}}.
  \end{align*}
Since $p_1$ does not depend on $k$, $T_0$ can also be chosen to be independent of $k\geq1$.
Substituting the estimate \eqref{sect-4-existence.2} into the above inequality and
by an analogous argument of \eqref{sect-4-existence.2.5}, we can find two constants
$C'_{1,p},C'_{2,p}>0$ and $T_1\leq T_0$, still
independent on $k$, such that for all $t\leq T_1$,
  \begin{align}\label{sect-4-existence.3}
  \hskip-2pt\|\rho^k_t\|_{L^p(\P\times\mu)}&\leq
  C'_{1,p}\bigg[\int_{\R^{n_1}} \exp\big\{C'_{2,p}T_1\big([\div_{x_1}(b_1)]^-
  +|\bar b_1|+|\bar\sigma_1|^2+|\nabla_{x_1}\sigma_1|^2\big)\big\}\d\mu_1\bigg]^{\frac1{2p(p+1)}}\cr
  \hskip-2pt&\hskip12pt \times\bigg[\int_{\R^{n}} \exp\big\{ C'_{2,p}T_1\big([\div_{x_2}(b_2)]^-
  +|\bar b_2|+|\bar\sigma_2|^2+|\nabla_{x_2}\sigma_2|^2\big)\big\}\d\mu \bigg]^{\frac1{2p(p+1)}}.
  \end{align}

{\it Step 2}. We show in this step that the family of flows $(X^k_t)_{k\geq 1}$
are convergent in some sense. For this purpose we check the conditions of
Theorem \ref{sect-4-stability}. First, by Remark \ref{sect-2-rem}(ii), the inequality
\eqref{sect-2-prop-1.2} shows that \eqref{sect-4-stability.1} is satisfied.
Next by \eqref{sect-4-existence.3}, we see that under the assumptions (H2) and (H4),
 \begin{equation}\label{sect-4-existence.4}
 \Lambda_{p,T_1}:=\sup_{k\geq1}\sup_{0\leq t\leq T_1}\|\rho^k_t\|_{L^p(\P\times\mu)}<+\infty,
 \end{equation}
which is nothing but \eqref{sect-4-stability.3}. It remains to check \eqref{sect-4-stability.2}.
Similar to the proof of \eqref{sect-2-lem-3.1.5}, we have
  $$|\nabla_{x_2}b_{2,k}|\leq |\nabla_{x_2}b_{2}|\ast\chi_k+2C|\bar b_2|\ast\chi_k.$$
Thus
  \begin{equation}\label{sect-4-existence.5}
  \int_{B(R)}|\nabla_{x_2}b_{2,k}|^q\,\d x\leq
  C_q\int_{B(R)}\big[(|\nabla_{x_2}b_2|\ast\chi_k)^q+(|\bar b_2|\ast\chi_k)^q\big]\,\d x.
  \end{equation}
By Jensen's inequality,
  \begin{align*}
  \int_{B(R)}|\nabla_{x_2}b_{2,k}|^q\,\d x
  &\leq C_q\int_{B(R)}\big(|\nabla_{x_2}b_2|^q+|\bar b_2|^q\big)\ast\chi_k\,\d x\cr
  &\leq C_q \|\nabla_{x_2}b_2\|_{L^q(B(R+1))}^q+ C_q \|\bar b_2\|_{L^q(B(R+1))}^q.
  \end{align*}
Therefore
  $$\sup_{k\geq1}\|\nabla_{x_2}b_{2,k}\|_{L^q(B(R))}<+\infty.$$
Analogously, we can show that $\sup_{k\geq1}\|\nabla_{x_2}\sigma_{2,k}\|_{L^{2q}(B(R))}<+\infty$.
Hence \eqref{sect-4-stability.2} is also satisfied. By Theorem \ref{sect-4-stability}, there exists
$X_2:\Omega\times\R^n\ra C([0,T_1],\R^{n_2})$ such that
  \begin{equation}\label{sect-4-existence.6}
  \lim_{k\ra\infty}\E\int_{\R^n}1\wedge \|X^k_{2,\cdot}-X_{2,\cdot}\|_{\infty,T_1}\,\d\mu=0.
  \end{equation}

{\it Step 3}. In the last step we prove that the random field $X_t=(X_{1,t},X_{2,t})$ is
the stochastic flow generated by the It\^o SDE \eqref{Ito-SDE}. First the same proof
as that of Proposition \ref{sect-2-prop-2} shows that there exists a family
$\{\rho_t:0\leq t\leq T_1\}$ of density functions such that $(X_t)_\#\mu=\rho_t\mu$
for any $t\in[0,T_1]$. Moreover $\sup_{0\leq t\leq T_1}\|\rho_t\|_{L^p(\P\times\mu)}
\leq \Lambda_{p,T_1}$, where $\Lambda_{p,T_1}$ is defined in \eqref{sect-4-existence.4}.

Thanks to \eqref{sect-4-existence.6}, we have the following analogues of
Corollary \ref{sect-2-cor-2}:
  $$\lim_{k\ra\infty}\int_{\R^n}\E\bigg(\sup_{0\leq t\leq T}\bigg|\int_0^t
  \big[\sigma_{2,k}(X^k_s)-\sigma_2(X_s)\big]\d B_s\bigg|\bigg)\d\mu=0$$
and
  $$\lim_{k\ra\infty}\int_{\R^n}\E\bigg(\sup_{0\leq t\leq T}\bigg|\int_0^t
  \big[b_{2,k}(X^k_s)-b_2(X_s)\big]\d s\bigg|\bigg)\d\mu=0.$$
With the above two limit results in hand, we let $k$ goes to $+\infty$ in the following
equation
  $$X^k_{2,t}=x_2+\int_0^t\sigma_{2,k}(X^k_s)\,\d B_s+\int_0^tb_{2,k}(X^k_s)\,\d s$$
and conclude that $X_t$ is the flow generated by \eqref{Ito-SDE}.  \fin

\medskip

Now we show the uniqueness of generalized stochastic flow associated to It\^o SDE
\eqref{Ito-SDE} on the time interval $[0,T_1]$.

\begin{proposition}[Uniqueness]\label{sect-4-uniqueness}
Under the assumptions (H1), (H2), (H3$'$) and (H4), there is at most one generalized
stochastic flow associated to the It\^o SDE \eqref{Ito-SDE} on the interval $[0,T_1]$.
\end{proposition}

\noindent{\bf Proof.} Suppose there are two flows $X_t=(X_{1,t},X_{2,t})$ and
$\hat X_t=(X_{1,t},\hat X_{2,t})$ associated to \eqref{Ito-SDE}, such that
$(X_t)_\#\mu=\rho_t\mu$ and $(\hat X_t)_\#\mu=\hat\rho_t\mu$. Let
  $$\Lambda_{p,T_1}:=\sup_{0\leq t\leq T_1}\|\rho_t\|_{L^p(\P\times\mu)}
  \vee\|\hat\rho_t\|_{L^p(\P\times\mu)}$$
which is finite. Applying Lemma \ref{a-priori-estimate}, we have
  \begin{align}\label{uniqueness.1}
  &\E\int_{G_R\cap\hat G_R}\log\bigg(\frac{\|X_{2,\cdot}-\hat X_{2,\cdot}\|^2_{\infty,T_1}}{\delta^2}+1\bigg)\d\mu\cr
  &\hskip6pt \leq C_{T_1}\Lambda_{p,T_1}C_{n_2,q}\Big[\|\nabla_{x_2}b_2\|_{L^q(B(4R))}
  +\|\nabla_{x_2}\sigma_2\|_{L^{2q}(B(4R))}+\|\nabla_{x_2}\sigma_2\|^2_{L^{2q}(B(4R))}\Big],
  \end{align}
where $G_R$ (resp. $\hat G_R$) is the level set of $X_t$ (resp. $\hat X_t$). Fix $R>0$, we see that
the right hand side is bounded, independent of $\delta>0$. Define, for $\eta>0$,
  $$\Sigma_\eta=\big\{(\omega,x)\in\Omega\times\R^n:
  \|X_{2,\cdot}(\omega,x)-\hat X_{2,\cdot}(\omega,x)\|_{\infty,T_1}\geq\eta\big\}.$$
Then by \eqref{uniqueness.1}, we have
  \begin{align*}
  \E\int_{G_R\cap\hat G_R}\ch_{\Sigma_\eta}\,\d\mu
  &\leq \frac1{\log\big(\frac{\eta^2}{\delta^2}+1\big)}
  \E\int_{G_R\cap\hat G_R}
  \log\bigg(\frac{\|X_{2,\cdot}-\hat X_{2,\cdot}\|_{\infty,T_1}}{\delta^2}+1\bigg) \d\mu\cr
  &\leq \frac{\bar C_{n_2,q,R,T_1}}{\log\big(\frac{\eta^2}{\delta^2}+1\big)}.
  \end{align*}
Note that the right hand side goes to 0 as $\delta\da0$, hence
  $$\E\int_{G_R\cap\hat G_R}\ch_{\Sigma_\eta}\,\d\mu=0$$
for any fixed $\eta>0$. Let $\eta\da0$, we obtain
  \begin{equation}\label{uniqueness.2}
  (\P\times\mu)\big\{(\omega,x)\in G_R\cap\hat G_R:\|X_{2,\cdot}-\hat X_{2,\cdot}\|_{\infty,T_1}>0\big\}=0.
  \end{equation}

Now notice that under the hypotheses (H2) and (H4), the estimates of level sets in
Lemma \ref{sect-3-lem-1} still hold. Therefore
  $$(\P\times\mu)\big[(G_R\cap\hat G_R)^c\big]\leq (\P\times\mu)(G_R^c)
  +(\P\times\mu)(\hat G_R^c)\leq \frac CR.$$
From this inequality it is clear that $G_R\cap\hat G_R \ua \Omega\times\R^n$ as $R$ increases to $+\infty$.
Letting $R\ua+\infty$ in \eqref{uniqueness.2}, we see that $(\P\times\mu)$ a.s.,
$\|X_{2,\cdot}-\hat X_{2,\cdot}\|_{\infty,T_1}=0$. \fin

\medskip

Following the arguments of Section 2, we can finally extend the flow $X_t$ to
any time interval $[0,T]$; moreover, the push-forward $(X_t)_\#\mu=\rho_t\mu$ and
the density function $\rho_t\in L^1\log L^1$.

\section{Weak differentiability of generalized stochastic flow}

Using the results of the preceding section, we intend to prove
in this section that the generalized stochastic flow
associated to the It\^o SDE with Sobolev coefficients, for which the
existence and uniqueness were established in Theorem \ref{Extension}
(see also \cite{Zhang10, FangLuoThalmaier, Zhang12}),
is weakly differentiable in the sense of measure, as in \cite{LeBrisLions04}.

First we introduce some notations and assumptions. Let $d,m\geq1$ be integers. Suppose we are given
a matrix-valued function $\sigma:\R^d\ra\R^m\otimes\R^d$ and a vector field $b:\R^d\ra\R^d$.
$B_t$ is an $m$-dimensional standard Brownian motion. We consider the following
It\^o's SDE
  \begin{equation}\label{sect-5.1}
  \d X_t(x)=\sigma(X_t(x))\,\d B_t+b(X_t(x))\,\d t,\quad X_0(x)=x.
  \end{equation}
In this section we write $X_t(x)$ to stress the initial condition of the stochastic flow.
Fix $q>1$ and $\alpha_1>d/2$. We denote by $\d\mu_1(x)=(1+|x|^2)^{-\alpha_1}\,\d x$
which is a finite measure on $\R^d$. We still write $\bar \sigma$ (resp. $\bar b$)
for $\frac\sigma{1+|x|}$ \big(resp. $\frac b{1+|x|}$\big).
Our assumptions in this section are:
\begin{itemize}
\item[\rm(A1)] $\sigma\in W^{1,2q}_{loc}$ and $b\in W^{1,q}_{loc}$;
\item[\rm(A2)] $\int_{\R^d}\exp\big[p_0\big([\div(b)]^- +|\bar b|+|\bar\sigma|^2
+|\nabla\sigma|^2\big)\big]\,\d \mu_1<+\infty$  for some $p_0>0$.
\end{itemize}

By Theorem \ref{Extension}, we see that under the assumptions (A1) and (A2),
the SDE \eqref{sect-5.1} generates a unique stochastic flow $X_t$ of measurable maps
on $\R^d$, such that the reference measure $\mu_1$ is absolutely continuous
under the flow. In order to prove the weak differentiability of the map $X_t:\R^d\ra\R^d$,
we need one more condition:
\begin{itemize}
\item[\rm(A3)] $\int_{\R^d}e^{p_0|\nabla b|}\,\d \mu_1<+\infty$ for some $p_0>0$.
\end{itemize}

We follow the line of arguments in \cite[Section 4]{LeBrisLions04}. Consider the It\^o SDE
on $\R^{2d}$:
  \begin{equation}\label{sect-5.2}
  \begin{cases}
  \d X_t(x)=\sigma(X_t(x))\,\d B_t+b(X_t(x))\,\d t, & X_0(x)=x,\\
  \d Y_t(x,y)=\big[\nabla\sigma(X_t(x))\big]Y_t(x,y)\,\d B_t
  +\big[\nabla b(X_t(x))\big]Y_t(x,y)\,\d t, & Y_0(x,y)=y.
  \end{cases}
  \end{equation}
As mentioned for the case of ODE in \cite[Section 4]{LeBrisLions04},
the above system of equations should be the limit of a system obtained by
perturbing the initial condition of the first equation. That is, for $\ee>0$,
we may consider
  $$\d X_t(x+\ee y)=\sigma(X_t(x+\ee y))\,\d B_t+b(X_t(x+\ee y))\,\d t,\quad X_0(x+\ee y)=x+\ee y.$$
Combining this equation together with \eqref{sect-5.1}, we obtain a system:
  \begin{equation}\label{sect-5.3}
  \begin{cases}
  \d X_t(x)=\sigma(X_t(x))\,\d B_t+b(X_t(x))\,\d t, & X_0(x)=x,\\
  \d \big[\frac{X_t(x+\ee y)-X_t(x)}\ee\big]=\frac{\sigma(X_t(x+\ee y))-\sigma(X_t(x))}\ee\,\d B_t
  +\frac{ b(X_t(x+\ee y))- b(X_t(x))}\ee\,\d t, & \frac{X_0(x+\ee y)-X_0(x)}\ee=y.
  \end{cases}
  \end{equation}
Now it is clear that the system of equations \eqref{sect-5.2} should be the limit in a certain
sense of the above system as $\ee\ra0$.

We now interpret both systems \eqref{sect-5.2} and \eqref{sect-5.3} as
the It\^o SDE with partially Sobolev coefficients studied in Section 4:
  $$\begin{cases}
  \d X_{1,t}=\sigma_1(X_{1,t})\,\d B_t+b_1(X_{1,t})\,\d t, & X_{1,0}=x_1,\\
  \d X_{2,t}=\sigma_2(X_{1,t},X_{2,t})\,\d B_t+b_2(X_{1,t},X_{2,t})\,\d t, & X_{2,0}=x_2.
  \end{cases}$$
where $x=(x_1,x_2)\in\R^{n_1}\times\R^{n_2}$ and $n_1+n_2=n$. In fact,
\begin{itemize}
\item for system \eqref{sect-5.2}, we set $x_1=x,x_2=y, n_1=n_2=d, X_{1,t}=X_t,X_{2,t}=(\nabla_x X_t)\,y,
\sigma_1=\sigma,b_1=b$ and $\sigma_2=(\nabla_x\sigma)\, y, b_2=(\nabla_x b)\, y$;
\item for system \eqref{sect-5.3}, we introduce the parameter $\ee>0$ and set
$x_1=x,x_2=y, n_1=n_2=d, X_{1,t}=X_t,X^\ee_{2,t}=
\frac{X_t(x+\ee y)-X_t(x)}\ee, \sigma_1=\sigma,b_1=b$ and $\sigma^\ee_2=
\frac{\sigma(x+\ee y)-\sigma(x)}\ee,b^\ee_2=\frac{b(x+\ee y)-b(x)}\ee$.
\end{itemize}
In the following we shall show that the two systems \eqref{sect-5.2} and \eqref{sect-5.3}
interpreted as above verify the main conditions of Section 4, and that
the stochastic flows associated to \eqref{sect-5.3} are convergent to that of
\eqref{sect-5.2} as $\ee\ra0$. To this end, we shall fix $\alpha>2\alpha_1+q+d/2$ throughout
this section. The reason for this special choice of $\alpha$ will become clear in the
following proofs. Denote by
  $$\d\mu(x_1,x_2)=\frac{\d x_1\d x_2}{(1+|x_1|^2+|x_2|^2)^\alpha}.$$
Then $\mu$ is obviously a finite measure on $\R^{2d}$. We first prove

\begin{lemma}\label{sect-5-lem-1}
Under the assumptions (A1)--(A3), both systems \eqref{sect-5.2} and \eqref{sect-5.3}
satisfy the conditions (H1), (H2), (H3$'$) and (H4).
\end{lemma}

\noindent{\bf Proof.} First, note that for both systems \eqref{sect-5.2} and \eqref{sect-5.3},
the conditions (H1) and (H2) on $\sigma_1$ and $b_1$ are exactly the same
assumptions (A1) and (A2) for $\sigma$  and $b$. In the following
we check the hypotheses (H3$'$) and (H4) for the two systems under the additional
assumption (A3) on the drift vector field $b$.

(1) We first treat the system \eqref{sect-5.2}. Since $\sigma_2(x_1,x_2)=(\nabla\sigma(x_1))\,x_2$,
we have $\nabla_{x_2}\sigma_2(x_1,x_2)=\nabla\sigma(x_1)$, hence for any $R>0$,
  \begin{align*}
  &\int_{B_1(R)}\d x_1\int_{B_2(R)}\big(|\sigma_2(x_1,x_2)|^{2q}
  +|\nabla_{x_2}\sigma_2(x_1,x_2)|^{2q}\big)\,\d x_2\cr
  &\quad\leq \int_{B_1(R)}\d x_1\int_{B_2(R)}\big(|\nabla\sigma(x_1)|^{2q}|x_2|^{2q}
  +|\nabla\sigma(x_1)|^{2q}\big)\,\d x_2\cr
  &\quad \leq (1+R^{2q})\,\Sigma_d R^d \int_{B_1(R)}|\nabla\sigma(x_1)|^{2q}\,\d x_1<+\infty.
  \end{align*}
Recall that $B_i(R)$ is a ball in $\R^{n_i}=R^d$ centered at the origin with radius $R\
(i=1,2)$, and $\Sigma_d$ is the volume of unit ball in $\R^d$.
Hence $\sigma_2\in L^{2q}_{x_1,loc}\big(W^{1,2q}_{x_2,loc}\big)$. In the same way
we can show that $b_2\in L^{q}_{x_1,loc}\big(W^{1,q}_{x_2,loc}\big)$. As a result,
(H3$'$) is satisfied.

Next note that $\div_{x_2}(b_2)(x_1,x_2)=\div(b)(x_1)$ which is independent on $x_2\in\R^{n_2}=\R^d$.
Since $b_2(x_1,x_2)=(\nabla b(x_1))\,x_2$, we have
  $$|\bar b_2(x_1,x_2)|=\frac{|(\nabla b(x_1))\, x_2|}{1+|(x_1,x_2)|}\leq |\nabla b(x_1)|;$$
similarly $|\bar\sigma_2(x_1,x_2)|^2\leq |\nabla\sigma(x_1)|^2$.
Moreover, $|\nabla_{x_2}\sigma_2(x_1,x_2)|^2=|\nabla\sigma(x_1)|^2$.
Combining these facts, it is clear that the assumptions (A2) and (A3) imply that
$\sigma_2$ and $b_2$ satisfy the condition (H4) for some $p_1\in(0,p_0]$.

(2) Now we deal with the second system \eqref{sect-5.3}. First we show that
$b^\ee_2\in L^{q}_{x_1,loc}\big(W^{1,q}_{x_2,loc}\big)$ for any $\ee\leq 1$.
By Fubini's theorem,
  \begin{equation}\label{sect-5-lem-1.1}
  \int_{B_1(R)}\d x_1\int_{B_2(R)}|b^\ee_2(x_1,x_2)|^q\,\d x_2
  =\int_{B_2(R)}\d x_2\int_{B_1(R)}\ee^{-q}|b(x_1+\ee x_2)-b(x_1)|^q\,\d x_1.
  \end{equation}
For any fixed $\ee\leq 1$ and $x_2\in B_2(R)$, by the pointwise characterization
of Sobolev functions, we have for a.e. $x_1\in\R^{n_1}$,
  \begin{equation}\label{sect-5-lem-1.2}
  |b(x_1+\ee x_2)-b(x_1)|\leq C_d\,\ee|x_2|\big(M_{|x_2|}|\nabla b|(x_1+\ee x_2)
  +M_{|x_2|}|\nabla b|(x_1)\big).
  \end{equation}
Therefore
  \begin{align*}
  &\int_{B_1(R)}\ee^{-q}|b(x_1+\ee x_2)-b(x_1)|^q\,\d x_1\cr
  &\quad \leq C_{d,q}|x_2|^q \int_{B_1(R)}\big[\big(M_{|x_2|}|\nabla b|(x_1+\ee x_2)\big)^q
  +\big(M_{|x_2|}|\nabla b|(x_1)\big)^q\big]\,\d x_1.
  \end{align*}
For $\ee\leq 1$ and $|x_2|\leq R$, by the maximal function inequality,
  \begin{align*}
  \begin{split}
  \int_{B_1(R)}\big(M_{|x_2|}|\nabla b|(x_1+\ee x_2)\big)^q\,\d x_1
  &= \int_{\ee x_2+B_1(R)}\big(M_{|x_2|}|\nabla b|(u)\big)^q\,\d u\cr
  &\leq \int_{B_1(2R)}\big(M_R|\nabla b|(u)\big)^q\,\d u\cr
  &\leq C'_{d,q}\int_{B_1(3R)}|\nabla b(u)|^q\,\d u.
  \end{split}
  \end{align*}
Consequently,
  $$\int_{B_1(R)}\ee^{-q}|b(x_1+\ee x_2)-b(x_1)|^q\,\d x_1
  \leq \bar C_{d,q}|x_2|^q\int_{B_1(3R)}|\nabla b(u)|^q\,\d u.$$
Substituting this inequality into \eqref{sect-5-lem-1.1}, we easily see that
  $$\sup_{0<\ee\leq 1}\int_{B_1(R)}\d x_1\int_{B_2(R)}|b^\ee_2(x_1,x_2)|^q\,\d x_2
  \leq \bar C_{d,q}\Sigma_d R^{d+q}\|\nabla b\|^q_{L^q(B_1(3R))}<+\infty,$$
where $\Sigma_d$ is the volume of the unit ball in $\R^d$. Therefore,
$b^\ee_2\in L^q_{x_1,loc}\big(L^q_{x_2,loc}\big)$. Next, since
$\nabla_{x_2}b^\ee_2(x_1,x_2)=\nabla b(x_1+\ee x_2)$, it is easy to show that
$\nabla_{x_2}b^\ee_2\in L^q_{x_1,loc}\big(L^q_{x_2,loc}\big)$. Hence the assertion
follows. In the same way we can show that $\sigma^\ee_2\in L^{2q}_{x_1,loc}
\big(W^{1,2q}_{x_2,loc}\big)$ for any $\ee\leq 1$. Thus we have finished
verifying (H3$'$).

The verifications of (H4) for $\sigma^\ee_2$ and $b^\ee_2$ are more complicated.
First we have $\div_{x_2}(b^\ee_2)(x_1,x_2)=\div(b)(x_1+\ee x_2)$. Hence
for $p>0$,
  $$K_{1,\ee}:=\int_{\R^{2d}}e^{p[\div_{x_2}(b^\ee_2)]^-}\,\d\mu(x_1,x_2)
  =\int_{\R^d}\d x_2\int_{\R^d}\frac{e^{p[\div(b)(x_1+\ee x_2)]^-}}{(1+|x_1|^2+|x_2|^2)^\alpha}\,\d x_1.$$
Making the change of variable $u_1=x_1+\ee x_2$ in the inner integral leads to
  $$K_{1,\ee}=\int_{\R^d}\d x_2\int_{\R^d}\frac{e^{p[\div(b)(u_1)]^-}}{(1+|u_1-\ee x_2|^2+|x_2|^2)^\alpha}\,\d u_1.$$
When $\ee\leq 1/2$, one has $|u_1|^2\leq 2|u_1-\ee x_2|^2+2|\ee x_2|^2
\leq 2|u_1-\ee x_2|^2+|x_2|^2/2$, thus
  \begin{equation}\label{sect-5-lem-1.3}
  1+|u_1-\ee x_2|^2+|x_2|^2\geq (1+|u_1|^2+|x_2|^2)/2.
  \end{equation}
Therefore
  $$K_{1,\ee}\leq 2^\alpha\int_{\R^d}\d x_2\int_{\R^d}\frac{e^{p[\div(b)(u_1)]^-}}{(1+|u_1|^2+|x_2|^2)^\alpha}\,\d u_1
  \leq 2^\alpha \mu_2(\R^d)\int_{\R^d}e^{p[\div(b)(u_1)]^-}\,\d\mu_1(u_1),$$
where $\d\mu_2=(1+|x_2|^2)^{\alpha_1-\alpha}\,\d x_2$ is a finite measure
on $\R^{n_2}=\R^d$. Therefore by (A2), if $p\leq p_0$, we have
  \begin{equation}\label{sect-5-lem-1.4}
  \sup_{\ee\leq 1/2}\int_{\R^{2d}}e^{p[\div_{x_2}(b^\ee_2)]^-}\,\d\mu <+\infty.
  \end{equation}

We now prove that $\int_{\R^{2d}}e^{p|\bar b^\ee_2|}\,\d\mu<+\infty$ for
$p$ sufficiently small. In fact,
  $$\int_{\R^{2d}}e^{p|\bar b^\ee_2|}\,\d\mu=
  \int_{\R^d}\d x_2\int_{\R^d}\frac{\exp\Big\{p\frac{|b(x_1+\ee x_2)-b(x_1)|}{\ee(1+|(x_1,x_2)|)}\Big\}}
  {(1+|x_1|^2+|x_2|^2)^\alpha}\,\d x_1.$$
Again by the pointwise inequality \eqref{sect-5-lem-1.2}, we get
  \begin{equation}\label{sect-5-lem-1.5}
  \int_{\R^{2d}}e^{p|\bar b^\ee_2|}\,\d\mu
  \leq \int_{\R^d}\d x_2\int_{\R^d}\frac{\exp\big\{pC_d\big(M_{|x_2|}|\nabla b|(x_1+\ee x_2)
  +M_{|x_2|}|\nabla b|(x_1)\big)\big\}}{(1+|x_1|^2+|x_2|^2)^\alpha}\,\d x_1.
  \end{equation}
We first estimate the term
  $$K_{2,\ee}:=\int_{\R^d}\d x_2\int_{\R^d}\frac{\exp\big\{pC_d M_{|x_2|}|\nabla b|(x_1+\ee x_2)
  \big\}}{(1+|x_1|^2+|x_2|^2)^\alpha}\,\d x_1.$$
Similar to the treatment of $K_{1,\ee}$, changing the variable and by
\eqref{sect-5-lem-1.3}, we have for all $\ee\leq 1/2$,
  \begin{align*}
  K_{2,\ee}&\leq 2^\alpha\int_{\R^d}\d x_2\int_{\R^d}\frac{\exp\big\{pC_d M_{|x_2|}|\nabla b|(u_1)
  \big\}}{(1+|u_1|^2+|x_2|^2)^\alpha}\,\d u_1\cr
  &\leq 2^\alpha\int_{\R^d}\frac{\d x_2}{(1+|x_2|^2)^{\alpha-\alpha_1}}
  \int_{\R^d}e^{pC_d M_{|x_2|}|\nabla b|(u_1)}\,\d\mu_1(u_1),
  \end{align*}
where the measure $\mu_1$ is defined at the beginning of this section.
We split the right hand side into two parts:
  \begin{align}\label{sect-5-lem-1.6}
  K_{2,\ee}&\leq 2^\alpha\int_{\{|x_2|\leq 1\}}\frac{\d x_2}{(1+|x_2|^2)^{\alpha-\alpha_1}}
  \int_{\R^d}e^{pC_d M_{|x_2|}|\nabla b|(u_1)}\,\d\mu_1(u_1)\cr
  &\quad +2^\alpha\int_{\{|x_2|>1\}}\frac{\d x_2}{(1+|x_2|^2)^{\alpha-\alpha_1}}
  \int_{\R^d}e^{pC_d M_{|x_2|}|\nabla b|(u_1)}\,\d\mu_1(u_1).
  \end{align}
Denoting the two terms by $K_{2,\ee}^{(1)}$ and $K_{2,\ee}^{(2)}$ respectively.
Now we are going to apply Lemma \ref{appendix-lem-2}. In the present
case, $\lambda(z)=(1+|z|^2)^{-\alpha_1}\ (z\in\R^d)$ and $\delta=1$ or $|x_2|$.
It is easy to show that for any $\delta\geq 1$,
  $$\Lambda_0=\sup_{k\geq1}\bigg(\frac{1+(k+1)^2\delta^2}{1+(k-1)^2\delta^2}
  \bigg)^{\alpha_1}= (1+4\delta^2)^{\alpha_1}.$$
Thus for $|x_2|>1$, an application of \eqref{appendix-lem-2.2} gives us
  \begin{align}\label{sect-5-lem-1.7}
  \int_{\R^d}e^{pC_d M_{|x_2|}|\nabla b|}\,\d\mu_1
  &\leq \int_{\R^d}\big(1+pC_d M_{|x_2|}|\nabla b|\big)\,\d\mu_1\cr
  &\quad +6 \cdot 5^d(1+4|x_2|^2)^{\alpha_1}\int_{\R^d} e^{2pC_d|\nabla b|}\,\d\mu_1.
  \end{align}
By Cauchy's inequality and \eqref{appendix-lem-2.1}, we obtain
  \begin{align}\label{sect-5-lem-1.8}
  \int_{\R^d} M_{|x_2|}|\nabla b|\,\d\mu_1&\leq
  \bigg[\mu_1(\R^d)\int_{\R^d} \big(M_{|x_2|}|\nabla b|\big)^2\,\d\mu_1\bigg]^{\frac12}\cr
  &\leq \bigg[24\cdot 5^d\mu_1(\R^d)(1+4|x_2|^2)^{\alpha_1}\int_{\R^d} |\nabla b|^2\,\d\mu_1\bigg]^{\frac12}\cr
  &= C'_d\|\nabla b\|_{L^2(\mu_1)}(1+4|x_2|^2)^{\alpha_1/2}.
  \end{align}
Substituting \eqref{sect-5-lem-1.8} into \eqref{sect-5-lem-1.7}, we can find
some positive constant $C_{p,d}>0$ such that
  \begin{equation*}
  \int_{\R^d}e^{pC_d M_{|x_2|}|\nabla b|}\,\d\mu_1
  \leq C_{p,d}(1+4|x_2|^2)^{\alpha_1}\int_{\R^d} e^{2pC_d|\nabla b|}\,\d\mu_1.
  \end{equation*}
Therefore
  \begin{equation*}
  K_{2,\ee}^{(2)}\leq 2^\alpha C_{p,d}
  \bigg(\int_{\R^d} e^{2pC_d|\nabla b|}\,\d\mu_1\bigg)
  \int_{|x_2|>1}\frac{(1+4|x_2|^2)^{\alpha_1}}{(1+|x_2|^2)^{\alpha-\alpha_1}}\,\d x_2.
  \end{equation*}
Since $\alpha>2\alpha_1+ d/2$, the second integral is finite. As a result,
  \begin{equation*}
  \sup_{\ee\leq 1/2}K_{2,\ee}^{(2)}\leq 2^\alpha \tilde C_{p,d}
  \int_{\R^d} e^{2pC_d|\nabla b|}\,\d\mu_1.
  \end{equation*}
By (A3), we see that when $p\leq p_0/(2C_d)$, the right hand side is finite. Notice that
  $$K_{2,\ee}^{(1)}\leq 2^\alpha \Sigma_d\int_{\R^d}e^{pC_d M_1|\nabla b|(u_1)}\,\d\mu_1(u_1),$$
where $\Sigma_d$ is the volume of the $d$-dimensional unit ball.
In the same way we can prove that $\sup_{\ee\leq 1/2}K_{2,\ee}^{(1)}<+\infty$ for
$p\leq p_0/(2C_d)$. Substituting these estimates into \eqref{sect-5-lem-1.6},
we conclude that if $p\leq p_0/(2C_d)$, $K_{2,\ee}$ is bounded uniformly in
$\ee\leq 1/2$.
The same computations lead to
  $$\sup_{\ee\leq 1}\int_{\R^d}\d x_2\int_{\R^d}\frac{\exp\big\{pC_d M_{\ee|x_2|}|\nabla b|(x_1)
  \big\}}{(1+|x_1|^2+|x_2|^2)^\alpha}\,\d x_1<+\infty.$$
Therefore an application of Cauchy's inequality to \eqref{sect-5-lem-1.5}
gives us that for any $p\leq p_0/(4C_d)$,
  \begin{equation}\label{sect-5-lem-1.9}
  \sup_{\ee\leq 1/2}\int_{\R^{2d}}e^{p|\bar b^\ee_2|}\,\d\mu<+\infty.
  \end{equation}

Analogously, we can show that when $p$ is small enough, it holds
  \begin{equation}\label{sect-5-lem-1.10}
  \sup_{\ee\leq 1/2}\int_{\R^{2d}}e^{p|\bar \sigma^\ee_2|^2}\,\d\mu<+\infty.
  \end{equation}
Finally, since $\nabla_{x_2}\sigma^\ee_2(x_1,x_2)=(\nabla\sigma)(x_1+\ee x_2)$, we
follow the arguments for estimating $K_{1,\ee}$ and arrive at
  $$\sup_{\ee\leq 1/2}\int_{\R^{2d}}e^{p|\nabla_{x_2}\sigma^\ee_2|^2}\,\d\mu<+\infty$$
for $p$ sufficiently small. Combining this estimate with \eqref{sect-5-lem-1.4},
\eqref{sect-5-lem-1.9} and \eqref{sect-5-lem-1.10}, we conclude that
$\sigma^\ee_2$ and $b^\ee_2$ satisfy the condition (H4), uniformly in $\ee\in(0,1/2]$.  \fin

\medskip

By Lemma \ref{sect-5-lem-1}, we can apply the main results of Section 4 (Theorem
\ref{sect-4-existence} and Proposition \ref{sect-4-uniqueness}) to both systems
\eqref{sect-5.2} and \eqref{sect-5.3}. Therefore, the system \eqref{sect-5.2} (resp. \eqref{sect-5.3})
generates a unique stochastic flow $Z_t(x,y)=(X_t(x),Y_t(x,y))$ \big(resp. $Z^\ee_t(x,y)=
\big(X_t(x),\ee^{-1}(X_t(x+\ee y)-X_t(x))\big)$\big);
moreover the Radon--Nikodym densities $\rho_t=\frac{\d(Z_t)_\#\mu}{\d\mu}$
and $\rho^\ee_t=\frac{\d(Z^\ee_t)_\#\mu}{\d\mu}$ exist, and there is a
$T_0>0$ small enough (note that by the uniform estimate in Lemma \ref{sect-5-lem-1},
$T_0$ does not depend on $\ee\leq 1/2$) such that
  \begin{equation}\label{sect-5.4}
  \Lambda_{p,T_0}:=\bigg(\sup_{0\leq t\leq T_0}\|\rho_t\|_{L^p(\P\times\mu)}\bigg)\bigvee
  \bigg(\sup_{\ee\leq 1/2}\sup_{0\leq t\leq T_0}\|\rho^\ee_t\|_{L^p(\P\times\mu)}\bigg)
  <+\infty,
  \end{equation}
where $p$ is the conjugate number of $q$. Next we want to prove that
$Y^\ee_t(x,y):=\ee^{-1}(X_t(x+\ee y)-X_t(x))$
is convergent to $Y_t(x,y)$ in a certain sense,
following the idea of Theorem \ref{sect-4-stability}.

\begin{theorem}\label{sect-5-thm}
Under the assumptions (A1)--(A3), we have for any $T>0$,
  $$\lim_{\ee\ra0}\E\int_{\R^{2d}}1\wedge\|Y^\ee_\cdot- Y_\cdot\|_{\infty,T}\,\d\mu=0.$$
\end{theorem}

\noindent{\bf Proof.} First we show that
  \begin{equation}\label{sect-5-thm.1}
  \lim_{\ee\ra0}\E\int_{\R^{2d}}1\wedge\|Y^\ee_\cdot- Y_\cdot\|_{\infty,T_0}\,\d\mu=0.
  \end{equation}
The proof is similar to that of Theorem \ref{sect-4-stability}, and we shall apply
Lemma \ref{a-priori-estimate} to show the convergence. It is easy to see that
for any $R>0$,
  \begin{equation*}
  \|\nabla_{x_2}b_2\|_{L^q(B(R))}+\|\nabla_{x_2}\sigma_2\|_{L^{2q}(B(R))}<+\infty.
  \end{equation*}
Noticing that we already have the uniform density estimate \eqref{sect-5.4},
hence it only remains to check the following conditions:
  \begin{equation}\label{sect-5-thm.2}
  C_1:=\sup_{\ee\leq 1/2}\big(\|\sigma^\ee_2\|_{L^{2q}(\mu)}
  +\|b^\ee_2\|_{L^{2q}(\mu)}\big)<+\infty
  \end{equation}
and
  \begin{equation}\label{sect-5-thm.3}
  \sigma^\ee_2\ra \sigma_2\mbox{ in }L^{2q}_{loc}(\R^{2d}) \quad
  \mbox{and} \quad b^\ee_2\ra b_2\mbox{ in }L^{q}_{loc}(\R^{2d}).
  \end{equation}

By Remark \ref{sect-2-rem} and \eqref{sect-5-lem-1.9}, \eqref{sect-5-lem-1.10},
we easily deduce that $C_1$ defined in \eqref{sect-5-thm.2} is finite.
Next, since $\sigma^\ee_2(x_1,x_2)=\frac{\sigma(x_1+\ee x_2)-\sigma(x_1)}\ee$
and $\sigma_2(x_1,x_2)=(\nabla\sigma(x_1))\, x_2$, the convergence $\sigma^\ee_2\ra \sigma_2$
in $L^{2q}_{loc}(\R^{2d})$ follows from the fact that $\sigma\in W^{1,2q}_{loc}(\R^d)$.
Similarly we conclude that $b^\ee_2$ converge to $ b_2$ in $L^{q}_{loc}(\R^{2d})$.
Hence the convergences in \eqref{sect-5-thm.3} are verified. Now we are ready
to follow the line of the proof of Theorem \ref{sect-4-stability} to
obtain the convergence \eqref{sect-5-thm.1}.

We then follow the arguments of Proposition \ref{sect-2-prop-4} and
use the flow properties of $Z_t=(X_t,Y_t)$ and $Z^\ee_t=(X_t,Y^\ee_t)$ to
extend the convergence to the whole interval $[0,T]$. \fin

\medskip

This theorem shows that the generalized stochastic flow associated to the
It\^o SDE \eqref{sect-5.1} is weakly differentiable in the sense of measure,
provided that its coefficients $\sigma$ and $b$ satisfy the assumptions
(A1)--(A3). Note that if $\sigma$ and $b$ are globally Lipschitz continuous,
then they fulfil (A1)--(A3). In this case, however, our result is weaker
than that in \cite{BouleauHirsch}, where the authors proved that almost
surely, the map $X_t:\R^d\ra\R^d$ is almost everywhere differentiable with respect to the
initial data for any time, by using the theory of Dirichlet form. In \cite[Section 5]{LiLuo},
we considered the Stratonovich SDE with smooth diffusion coefficient $\sigma$
and Sobolev drift coefficient $b$, and proved the approximate differentiability
of the generalized stochastic flow by using the Ocone-Pardoux decomposition,
which essentially reduces the problem to prove the differentiability of the flow
generated by some ODE with random Sobolev coefficient.

\section{Appendix}

In this section we present some results that are used in the paper. We assume the
coefficients $\sigma:\R^n\ra\R^m\otimes\R^n$ and $b:\R^n\ra\R^n$ of the It\^o SDE
  \begin{equation}\label{appendix.1}
  \d X_t=\sigma(X_t)\,\d B_t+b(X_t)\,\d t,\quad X_0=x
  \end{equation}
are smooth and bounded together with their derivatives of all orders. Here $B_t$ is still an
$m$-dimensional standard Brownian motion. Then the above equation
generates a stochastic flow $X_t$ of diffeomorphisms on $\R^n$.

First we recall the
expression for the Radon--Nikodym density of the stochastic flow with respect to
some reference measure. Let $\lambda\in C^2(\R^n)$ and define a measure on $\R^n$ by
  $$\d\mu(x)=e^{\lambda(x)}\d x.$$
It is well known that the push-forward $(X_t)_\#\mu$
(resp. $(X_t^{-1})_\#\mu$) of $\mu$ by the flow $X_t$ (resp. the inverse flow $X_t^{-1}$) is
absolutely continuous with respect to $\mu$. Denote by
  $$\rho_t(x)=\frac{\d[(X_t)_\#\mu](x)}{\d\mu(x)} \quad \mbox{and} \quad
  \tilde\rho_t(x)=\frac{\d[(X_t^{-1})_\#\mu](x)}{\d\mu(x)}.$$
We have the following simple identity:
  \begin{equation}\label{appendix.1.5}
  \rho_t(x)=1/\tilde\rho_t\big(X_t^{-1}(x)\big).
  \end{equation}
Moreover by \cite[Lemma 4.3.1]{Kunita90}, a simple computation gives us (see also
\cite[(3.6)]{Zhang12})
  \begin{equation}\label{appendix.2}
  \tilde\rho_t(x)=\exp\bigg(\int_0^t\<\Lambda_1^\sigma(X_s(x)),\d B_s\>
  +\int_0^t\Lambda_2^{\sigma,b}(X_s(x))\,\d s\bigg),
  \end{equation}
in which
  $$\Lambda_1^\sigma=\div(\sigma)+\sigma^\ast\nabla\lambda \quad\mbox{and}\quad
  \Lambda_2^{\sigma,b}=\div(b)+\L\lambda-\frac12\<\nabla\sigma,(\nabla\sigma)^\ast\>.$$
Here by $\div(\sigma)=\big(\div(\sigma^{\cdot,1}),\ldots,\div(\sigma^{\cdot,m})\big)$
we mean the $\R^m$-valued function whose components are the divergences of the
columns of $\sigma$; $\sigma^\ast$ is the transpose of $\sigma$ and $\L$ is the
second order differential operator associated to \eqref{appendix.1}:
  $$\L\lambda=\frac12\sum_{i,j=1}^n a^{ij}\partial_i\partial_j\lambda
  +\sum_{i=1}^n b^i\partial_i\lambda$$
with $a^{ij}=\sum_{k=1}^m\sigma^{ik}\sigma^{jk}$ and $\partial_i\lambda=
\frac{\partial}{\partial x^i}\lambda$. Finally
  $$\<\nabla\sigma,(\nabla\sigma)^\ast\>=\sum_{k=1}^m\<\nabla\sigma^{\cdot,k},
  (\nabla\sigma^{\cdot,k})^\ast\>=\sum_{k=1}^m\sum_{i,j=1}^n(\partial_i\sigma^{jk})
  (\partial_j\sigma^{ik}).$$
From this expression, we see that if the first $n_1$-rows $\sigma_1=(\sigma^{ij})_{1\leq i\leq n_1,
1\leq j\leq m}$ only depend on the variables $x_1=(x^1,\ldots,x^{n_1})$, then
  \begin{align}\label{key-observation}
  \<\nabla\sigma,(\nabla\sigma)^\ast\>&=\sum_{k=1}^m\bigg(\sum_{i,j=1}^{n_1}(\partial_i\sigma^{jk})
  (\partial_j\sigma^{ik})+\sum_{i,j=n_1+1}^{n}(\partial_i\sigma^{jk})
  (\partial_j\sigma^{ik})\bigg)\cr
  &=\<\nabla_{x_1}\sigma_1,(\nabla_{x_1}\sigma_1)^\ast\>
  +\<\nabla_{x_2}\sigma_2,(\nabla_{x_2}\sigma_2)^\ast\>,
  \end{align}
where $\sigma_2$ consists of the last $(n-n_1)$-rows of the matrix $\sigma$.
Notice that the derivatives $\nabla_{x_1}\sigma_2$ are not involved here.
This observation is crucial for the present work.

The following is an $L^p$-estimate for $\rho_t(x)$ which is proved in \cite[Lemma 3.2]{Zhang12}
(see also \cite[Theorem 2.1]{FangLuoThalmaier} for the case where $\mu=\gamma_n$
is the standard Gaussian measure).

\begin{lemma}\label{density-estimate}
Assume that $\mu(\R^n)<+\infty$. Then for any $t\in [0,T]$ and $p>1$,
  \begin{equation}\label{density-estimate.1}
  \|\rho_t\|_{L^p(\P\times\mu)}\leq \mu(\R^n)^{1/(p+1)}
  \bigg(\sup_{t\in[0,T]}\int_{\R^n}\exp\big(tp^3|\Lambda_1^\sigma|^2
  -tp^2\Lambda_2^{\sigma,b}\big)\d\mu\bigg)^{1/p(p+1)}.
  \end{equation}
\end{lemma}

Next we present a simple technical result.

\begin{lemma}\label{appendix-lem} Let $f\in L_{loc}^1(\R^n)$ and denote by
$\bar f=\frac{f}{1+|x|}$. Then
  \begin{equation}\label{appendix-lem.1}
  \frac{|f\ast\chi_k|(x)}{1+|x|} \leq 2 (|\bar f|\ast\chi_k)(x),\quad x\in\R^n.
  \end{equation}
\end{lemma}

\noindent{\bf Proof.}
Indeed, for each $k\geq1$,
  $$\frac{|f\ast\chi_k|(x)}{1+|x|}\leq
  \int_{B(1/k)}\frac{|f(x-y)|}{1+|x|}\,\chi_k(y)\,\d y.$$
For $y\in B(1/k)$, one has $|x-y|\leq |x|+1/k$, thus
  $$1+|x-y|\leq 2+|x|\leq 2(1+|x|).$$
As a result,
  \begin{align*}
  \frac{|f\ast\chi_k|(x)}{1+|x|}
  \leq 2\int_{B(1/k)}\frac{|f(x-y)|}{1+|x-y|}\,\chi_k(y)\,\d y
  = 2\int_{B(1/k)}|\bar f(x-y)|\,\chi_k(y)\,\d y,
  \end{align*}
from which we deduce \eqref{appendix-lem.1}. \fin

\medskip

In the following we introduce the pointwise inequality for partially Sobolev functions.
To this end, we need the notion of locally maximal function for partial variables.
As in the introduction, $n=n_1+n_2$ and for $x\in\R^n$, we write $x=(x_1,x_2)$
where $x_1\in\R^{n_1}$ and $x_2\in\R^{n_2}$. Let $f:\R^{n_1}\times\R^{n_2}\ra \R$
be locally integrable. For almost every $x_1\in\R^{n_1}$, define
  \begin{align*}
  M_{2,R}f(x_1,x_2)&=\sup_{0<r\leq R}\bint_{B_2(x_2,r)}|f(x_1,y_2)|\,\d y_2\cr
  &:=\sup_{0<r\leq R}\frac1{\L_{n_2}(B_2(x_2,r))}\int_{B_2(x_2,r)}|f(x_1,y_2)|\,\d y_2,\quad R>0.
  \end{align*}
Here $B_2(x_2,r)$ means the ball in $\R^{n_2}$ centered at $x_2$ with radius $r$.
Recall that $B_i(r)$ is the ball in $\R^{n_i}$ of radius $r$ centered at the origin, $i=1,2$.
The main point of the first result in the next lemma lies in the fact that the exceptional set $N$ is chosen
to be a negligible subset of $\R^n$.

\begin{lemma}\label{pointwise-inequality} \quad
\begin{itemize}
\item[\rm(i)] Suppose that $f:\R^{n_1}\times\R^{n_2}\ra \R$ belongs to the space
$L^1_{x_1,loc}\big(W^{1,1}_{x_2,loc}\big)$. Then there is a dimensional constant $C>0$
(independent of $f$) and a negligible set $N\subset \R^{n_1}\times\R^{n_2}$,
such that for all $(x_1,x_2),(x_1,y_2)\notin N$ with $|x_2-y_2|_{\R^{n_2}}\leq R$,
it holds
  $$|f(x_1,x_2)-f(x_1,y_2)|\leq C|x_2-y_2|_{\R^{n_2}}
  \big[M_{2,R}|\nabla_{x_2}f|(x_1,x_2)+M_{2,R}|\nabla_{x_2}f|(x_1,y_2)\big].$$
\item[\rm(ii)] If $f\in L^p_{loc}(\R^{n_1}\times\R^{n_2})$ for some $p>1$, then there is a constant
$C_{p,n_2}>0$ such that
  $$\int_{B_2(r)}(M_{2,R}f(x_1,x_2))^p\,\d x_2
  \leq C_{p,n_2}\int_{B_2(r+R)}|f(x_1,x_2)|^p\,\d x_2.$$
\end{itemize}
\end{lemma}

\noindent{\bf Proof.} (i) Here we present a proof based on the well known
pointwise inequality for Sobolev functions. Let
  $$\tilde N=\bigg\{(x_1,x_2)\in\R^n:x_1\in\R^{n_1} \mbox{ and } \limsup_{\L_{n_2}(B)\ra0,\, x_2\in B}
  \bigg|\bint_B f(x_1,y_2)\,\d y_2 -f(x_1,x_2)\bigg|>0\bigg\},$$
where the limit is taken over all balls $B\subset\R^{n_2}$ such that $x_2$ is contained in $B$.
$\tilde N$ is a measurable subset of $\R^n$. We see that for all $x_1\in\R^{n_1}$, the section
  $$\tilde N_{x_1}=\bigg\{x_2\in\R^{n_2}:\limsup_{\L_{n_2}(B)\ra0,\, x_2\in B}
  \bigg|\bint_B f(x_1,y_2)\,\d y_2 -f(x_1,x_2)\bigg|>0\bigg\}.$$
Since $f\in L^1_{x_1,loc}\big(W^{1,1}_{x_2,loc}\big)$, there is an $\L_{n_1}$-negligible
set $N_1\subset\R^{n_1}$, such that for every
$x_1\notin N_1$, one has $f(x_1,\cdot)\in W^{1,1}_{x_2,loc}$. In particular,
$f(x_1,\cdot)\in L^1_{x_2,loc}$. Lebesgue's differentiation
theorem gives us $\L_{n_2}(\tilde N_{x_1})=0$ for all $x_1\notin N_1$.
By Fubini's theorem we have
  $$\L_n(\tilde N)=\int_{\R^{n_1}}\L_{n_2}(\tilde N_{x_1})\,\d x_1=0.$$

Define $N=\tilde N\cup (N_1\times \R^{n_2})$. We see that $\L_n(N)=0$.
Now fix any $(x_1,x_2),(x_1,y_2)\notin N$ with $|x_2-y_2|_{\R^{n_2}}\leq R$. Since
$x_1\notin N_1$, we have $f(x_1,\cdot)\in W^{1,1}_{x_2,loc}$.
By the pointwise inequality of Sobolev functions (see e.g. \cite[p.186]{Ambrosio11}
or \cite[Theorem A.1]{FangLuoThalmaier}),
there exist a constant $C_{n_2}>0$ such that
for all $u_2,v_2\notin \tilde N_{x_1}$ with $|u_2-v_2|_{\R^{n_2}}\leq R$, it holds
  $$|f(x_1,u_2)-f(x_1,v_2)|\leq C|u_2-v_2|_{\R^{n_2}}
  \big[M_{2,R}|\nabla_{x_2}f|(x_1,u_2)+M_{2,R}|\nabla_{x_2}f|(x_1,v_2)\big].$$
Now the result follows by noticing that $x_2,y_2\notin N_{x_1}$ and
$\tilde N_{x_1}\subset N_{x_1}$.

(ii) This is obvious from the properties of maximal functions. \fin

\medskip

The next result is similar to Lemma \ref{pointwise-inequality}(ii), but the integral
is taken with respect to some other reference measure. Perhaps such a result
already exists, but we are unaware of its reference. We present its proof
for the reader's convenience. Suppose we are given a continuous $\lambda\in C(\R^n,(0,+\infty))$
such that $\d\mu=\lambda\,\d x$ is a finite measure on $\R^n$. Fix $\delta>0$.
For every positive integer $k$, we denote by $R_k:=\{x\in\R^n:(k-1)\delta\leq |x|\leq k\delta\}$, that is, the
ring between the concentric spheres centered at the origin with radii $(k-1)\delta$
and $k\delta$, respectively. Set
  $$\ol\lambda_k=\sup_{x\in R_k}\lambda(x),\quad
  \ul{\lambda}_k=\inf_{x\in (R_k)_\delta}\lambda(x),$$
where $(R_k)_\delta$ is the $\delta$-neighborhood
of the ring $R_k$. We shall denote by
  $$\Lambda_0=\sup_{k\geq1}\frac{\ol\lambda_k}{\ul\lambda_k}.$$
Obviously $\Lambda_0\geq1$. If $\lambda(x)=\phi(|x|)$ and for some $\beta>1$, $\phi(s)\sim e^{-s^\beta}$
as $s\ra\infty$, then $\Lambda_0=+\infty$. Therefore the following result
does not hold for the standard Gaussian measure.

The local maximal function $M_\delta f(x)$ of
a locally integrable function $f\in L^1_{loc}$ is defined as usual:
  $$M_\delta f(x)=\sup_{0<r\leq \delta}\bint_{B(x,r)}|f(y)|\,\d y
  :=\sup_{0<r\leq \delta}\frac1{\L_n(B(x,r))}\int_{B(x,r)}|f(y)|\,\d y.$$

\begin{lemma}\label{appendix-lem-2}
Assume that $\Lambda_0<+\infty$ and denote by $C_p=5^n2^p p/(p-1)$ for $p>1$. Then
  \begin{equation}\label{appendix-lem-2.1}
  \int_{\R^n}(M_\delta f)^p\,\d\mu\leq 3C_p\Lambda_0\int_{\R^n}|f|^p\,\d\mu.
  \end{equation}
As a result, for any $\theta>0$,
  \begin{equation}\label{appendix-lem-2.2}
  \int_{\R^n}e^{\theta M_\delta f}\,\d\mu\leq \int_{\R^n}(1+\theta M_\delta f)\,\d\mu
  +6 \cdot 5^n\Lambda_0\int_{\R^n}e^{2\theta|f|}\,\d\mu.
  \end{equation}
\end{lemma}

\noindent{\bf Proof.} Note that
  \begin{equation}\label{appendix-lem-2.3}
  \int_{\R^n}(M_\delta f)^p\,\d\mu=\sum_{k=1}^\infty \int_{R_k}(M_\delta f)^p\,\d\mu
  \leq \sum_{k=1}^\infty \ol\lambda_k\int_{R_k}(M_\delta f)^p\,\d x.
  \end{equation}
Next we follow the idea of \cite[Chap. I, Section 1]{Stein} to show that
for any $p>1$,
  \begin{equation}\label{appendix-lem-2.4}
  \int_{R_k}(M_\delta f)^p\,\d x\leq C_p\int_{(R_k)_\delta}|f|^p\,\d x,
  \end{equation}
where $C_p=2^p 5^np/(p-1)$. Indeed, for any $s>0$, we define $R_k(s)=\{x\in R_k:M_\delta f(x)>s\}$
(note that $s\ra \L_n(R_k(s))$ is the distribution function of $M_\delta f$ when
restricted on $R_k$).
Then similar to the argument on \cite[pp. 6--7]{Stein}, we have
  \begin{equation}\label{appendix-lem-2.5}
  \L_n(R_k(s))\leq \frac{2\cdot 5^n}s\int_{(R_k)_\delta\cap\{|f|>s/2\}}|f(y)|\,\d y.
  \end{equation}
Next it is easy to show that
  $$\int_{R_k}(M_\delta f)^p\,\d x=p\int_0^\infty s^{p-1}\L_n(R_k(s))\,\d s.$$
Substituting \eqref{appendix-lem-2.5} into the above equality and changing
the order of integration, we finally get
  $$\int_{R_k}(M_\delta f)^p\,\d x\leq \frac{5^n 2^p p}{p-1}\int_{(R_k)_\delta}|f(y)|^p\,\d y.$$

Now by \eqref{appendix-lem-2.4} and the definition of $\ul\lambda_k$, we have
  $$\int_{R_k}(M_\delta f)^p\,\d x\leq \frac{C_p}{\ul\lambda_k}\int_{(R_k)_\delta}|f|^p\,\d\mu.$$
Substituting this inequality into \eqref{appendix-lem-2.3}, we obtain
  $$\int_{\R^n}(M_\delta f)^p\,\d\mu\leq C_p\sum_{k=1}^\infty
  \frac{\ol\lambda_k}{\ul\lambda_k}\int_{(R_k)_\delta}|f|^p\,\d\mu
  \leq 3C_p\Lambda_0\int_{\R^n}|f|^p\,\d\mu.$$

Finally, by expanding the exponential function, we have
  \begin{equation}\label{appendix-lem-2.6}
  \int_{\R^n}e^{\theta M_\delta f}\,\d\mu=\int_{\R^n}(1+\theta M_\delta f)\,\d\mu
  +\sum_{k=2}^\infty\frac{\theta^k}{k!}\int_{\R^n}(M_\delta f)^k\,\d\mu.
  \end{equation}
Applying the inequality proved above, we get, for any $k\geq2$,
  $$\int_{\R^n}(M_\delta f)^k\,\d\mu
  \leq 3\Lambda_0\frac{5^n 2^k k}{k-1}\int_{\R^n}|f|^k\,\d\mu
  \leq 3\cdot 5^n\Lambda_0 2^{k+1}\int_{\R^n}|f|^k\,\d\mu.$$
Therefore,
  $$\sum_{k=2}^\infty\frac{\theta^k}{k!}\int_{\R^n}(M_\delta f)^k\,\d\mu
  \leq 6\cdot 5^n\Lambda_0 \sum_{k=2}^\infty\frac{(2\theta)^k}{k!}\int_{\R^n}|f|^k\,\d\mu
  \leq 6\cdot 5^n\Lambda_0\int_{\R^n}e^{2\theta|f|}\,\d\mu. $$
The proof is completed by substituting this inequality into \eqref{appendix-lem-2.6}. \fin

\bigskip

\noindent{\bf Acknowledgement.} The author is grateful to the financial supports
of the National Natural Science Foundation of China (No. 11101407),
and the Key Laboratory of Random Complex Structures and Data Science,
Academy of Mathematics and Systems Science, Chinese Academy of Sciences
(No. 2008DP173182).

\bigskip

{\sc Dejun Luo}

{\sc Institute of Applied Mathematics}

{\sc Academy of Mathematics and Systems Science}

{\sc Chinese Academy of Sciences }

{\sc Beijing 100190 }

{\sc China}

{\sc E-mail:} luodj@amss.ac.cn


\begin{thebibliography}{a23}

\bibitem{Ambrosio04} L. Ambrosio, {\it Transport equation and Cauchy problem
for BV vector fields}, Invent. Math. 158 (2004), 227--260.

\bibitem{Ambrosio11} L. Ambrosio, {\it The flow associated to weakly differentiable
vector fields: recent results and open problems}. A. Bressan et al. (eds.),
``Nonlinear Conservation Laws and Applications'', pp. 181--193,
The IMA Volumes in Mathematics and its Applications 153, Springer, 2011.

\bibitem{AmbrosioFigalli09} L. Ambrosio and A. Figalli, {\it On flows associated
to Sobolev vector fields in Wiener space: an approach \`a la DiPerna--Lions}.
J. Funct. Anal. 256 (2009), no. 1, 179--214.

\bibitem{ALM} L. Ambrosio, M. Lecumberry and S. Maniglia, {\it Lipschitz
regularity and approximate differentiability of the DiPerna--Lions
flow}. Rend. Sem. Mat. Univ. Padova, 114 (2005), 29--50.

\bibitem{BouleauHirsch} N. Bouleau and F. Hirsch, {\it On the derivability,
with respect to the initial data, of the solution of a stochastic differential
equation with Lipschitz coefficients}. S\'eminaire Th\'eorie du Potentiel No. 9,
Springer, Berlin, 1989.

\bibitem{CiprianoCruzeiro05} F. Cipriano and A.B. Cruzeiro, {\it Flows
associated with irregular $\R^d$-vector fields}. J. Diff. Equations
210 (2005), 183--201.

\bibitem{CrippadeLellis} G. Crippa and C. de Lellis, {\it Estimates and regularity
results for the DiPerna--Lions flows}. J. Reine Angew. Math. 616 (2008), 15--46.

\bibitem{DiPernaLions89} R.J. Di Perna and P.L. Lions, {\it Ordinary differential equations,
transport theory and Sobolev spaces.} Invent. Math. 98 (1989),
511--547.

\bibitem{Dumas} H.S. Dumas, F. Golse and P. Lochak, {\it Multiphase averaging for
generalized flows on manifolds}. Ergodic Theory Dynam. Systems 14 (1994), no. 1, 53--67.

\bibitem{EvansGariepy} L.C. Evans and R.F. Gariepy, {\it Measure theory and fine properties of
functions}. Studies in Advanced Math., CRC Press, London, 1992.

\bibitem{FangLiLuo} Shizan Fang, Huaiqian Li and Dejun Luo, {\it Heat semi-group and generalized
flows on complete Riemannian manifolds}. Bull. Sci. math. 135 (2011), 565--600.

\bibitem{FangLuo07} Shizan Fang and Dejun Luo, {\it Flow of homeomorphisms
and stochastic transport equations}. Stoch. Anal. Appl. 25 (2007),
1079--1108.

\bibitem{FangLuo10} Shizan Fang and Dejun Luo, {\it Transport equations and quasi-invariant flows on
the Wiener space}. Bull. Sci. Math. 134 (2010), 295--328.

\bibitem{FangLuoThalmaier} Shizan Fang, Dejun Luo and Anton Thalmaier,
{\it Stochastic differential equations with coefficients in Sobolev
spaces}. J. Funct. Anal. 259 (2010), no. 5, 1129--1168.

\bibitem{Figalli} A. Figalli, {\it Existence and uniqueness of
martingale solutions for SDEs with rough or degenerate
coefficients}. J. Funct. Anal. 254 (2008), 109--153.

\bibitem{Kunita90} H. Kunita, {\it Stochastic Flows and Stochastic
Differentail Equations}. Cambridge University Press, 1990.

\bibitem{LeBrisLions04} C. Le Bris and P.L. Lions, {\it Renormalized
solutions of some transport equations with partially $W^{1,1}$
velocities and applications}. Ann. Mat. Pura Appl. 183 (2004),
97--130.

\bibitem{LeBrisLions08} C. Le Bris and P.L. Lions, {\it Existence and
uniqueness of solutions to Fokker-Planck type equations with
irregular coefficients}. Comm. Partial Differential Equations 33
(2008), 1272--1317.

\bibitem{LiLuo} Huaiqian Li and Dejun Luo, {\it Quasi-invariant flow generated
by Stratonovich SDE with BV drift coefficients.} Stoch. Anal. Appl. 30 (2012),
258--284.

\bibitem{Luo10} Dejun Luo, {\it Well-posedness of Fokker--Planck type equations on the
Wiener space.} Infin. Dimens. Anal. Quantum Probab. Relat. Top. 13
(2010), no. 2, 273--304.

\bibitem{Luo12} Dejun Luo, {\it Fokker--Planck type equations with Sobolev diffusion coefficients
and BV drift coefficients}. Acta Math. Sin. (Engl. Ser.), to appear.

\bibitem{Stein} E.M. Stein, {\it Singular integrals and differentiability
properties of functions}. Princeton University Press, Princeton, New Jersey, 1970.

\bibitem{Zhang05}  Xicheng Zhang, {\it Homeomorphic flows for multi-dimensional SDEs
with non-Lipschitz coefficients}. Stochastic Process. Appl. 115
(2005), no. 3, 435--448.

\bibitem{Zhang10} Xicheng Zhang, {\it Stochastic flows of SDEs with
irregular coefficients and stochastic transport equations}. Bull.
Sci. Math. 134 (2010), 340--378.

\bibitem{Zhang11} Xicheng Zhang, {\it Quasi-invariant stochastic flows of SDEs
with non-smooth drifts on compact manifolds}. Stochastic Process. Appl. 121 (2011),
no. 6, 1373--1388.

\bibitem{Zhang12} Xicheng Zhang, {\it Well-posedness and large deviation for
degenerate SDEs with Sobolev coefficients}. Rev. Mat. Iberoam., to appear.
\end{thebibliography}
\end{document}